\newcommand{\idty}{{\leavevmode{\rm 1\mkern -5.4mu I}}}
\newcommand{\Ir}{Z\!\!\!Z}
\newcommand{\Ibb}[1]{ {\rm I\ifmmode\mkern
            -3.6mu\else\kern -.2em\fi#1}}
\newcommand{\ibb}[1]{\leavevmode\hbox{\kern.3em\vrule
     height 1.2ex depth -.3ex width .2pt\kern-.3em\rm#1}}
\newcommand{\Cx}{{\ibb C}}
\newcommand{\be}{\begin{eqnarray}}
\newcommand{\ee}{\end{eqnarray}}
\newcommand{\ldel}{\Delta_{\Omega^1}}
\newcommand{\rdel}{_{\Omega^1}\Delta}
\newcommand{\un}[1]{{\underline{#1}}}
\renewcommand{\d}{\mbox{d}}
\documentstyle[12pt]{article}
\begin{document}
\renewcommand{\theequation} {\arabic{section}.\arabic{equation}}

\begin{tabbing}
\hspace*{12cm}\= GOET-TP 95/95 \\
              \> August 1995
\end{tabbing}
\vskip.5cm

\centerline{\huge \bf NONCOMMUTATIVE GEOMETRY}
\vskip.8cm
\centerline{\huge \bf OF FINITE GROUPS}
\vskip1cm

\begin{center}
  {\large \bf Klaus Bresser, \ Folkert M\"uller-Hoissen}
      \vskip.3cm
     Institut f\"ur Theoretische Physik,
     Bunsenstr. 9, D-37073 G\"ottingen
     \vskip.5cm
     {\large \bf Aristophanes Dimakis}
     \vskip.3cm
     Department of Mathematics, University of Crete, GR-71409
     Iraklion, Greece
      \vskip.5cm
     {\large \bf Andrzej Sitarz}
     \vskip.3cm
     Department of Theoretical Physics, Jagiellonian University \\
     Reymonta 4, 30-059 Krak{\'o}w, Poland
\end{center}

\vskip1.cm

\begin{abstract}
\noindent
A finite set can be supplied with a group structure which can then
be used to select (classes of) differential calculi on it via the
notions of left-, right- and bicovariance. A corresponding framework
has been developed by Woronowicz, more generally for Hopf algebras
including quantum groups.
A differential calculus is regarded as the most basic structure needed
for the introduction of further geometric notions like linear
connections and, moreover, for the formulation of field theories
and dynamics on finite sets. Associated with each bicovariant first
order differential calculus on a finite group is a braid operator which
plays an important role for the construction of distinguished geometric
structures. For a covariant calculus, there are notions of invariance
for linear connections and tensors. All these concepts are explored for
finite groups and illustrated with examples.
Some results are formulated more generally for arbitrary associative
(Hopf) algebras. In particular, the problem of extension of a
connection on a bimodule (over an associative algebra) to tensor
products is investigated, leading to the class of `extensible
connections'. It is shown that invariance properties of an extensible
connection on a bimodule over a Hopf algebra are carried over to the
extension. Furthermore, an invariance property of a connection is
also shared by a `dual connection' which exists on the dual bimodule
(as defined in this work).
\end{abstract}

\newpage

\section{Introduction}
\setcounter{equation}{0}
Noncommutative geometry (see \cite{Conn94}, for example) replaces the
familiar arena of classical physics, a manifold supplied with
differential geometric structures, by an associative algebra $\cal A$
and algebraic structures on it. According to our point of view, the most
basic geometric structure in the framework of noncommutative geometry is
a `differential calculus' on $\cal A$ (see also \cite{Herm77}). It
allows the introduction of further geometric notions like linear
connections and, moreover, the formulation of field theories and
dynamics on finite sets.
\vskip.2cm

Though noncommutative geometry is designed to handle noncommutative
algebras, nontrivial structures already arise on commutative algebras
with non-standard differential calculi (see \cite{BDMH95} and
references given there). A commutative algebra of particular interest
in this context is the algebra of functions on a finite (or discrete)
set. A differential calculus on a finite set provides the latter with
a structure which may be viewed as a discrete counterpart to that of
a (continuous) differentiable manifold \cite{DMH94-graphs,DMHV}.
It has been shown in \cite{DMH94-graphs} that (first order)
differential calculi on discrete sets are in correspondence with
(di)graphs with at most two (antiparallel) arrows between any two
vertices. This relation with graphs and networks suggests applications
of the formalism to dynamics on networks \cite{DMHV} and the universal
dynamics considered in \cite{Mack95}, for example.
\vskip.2cm

A finite set can always be supplied with a group structure. The
left- and right-action of the group on itself can then be used to
distinguish certain differential calculi and geometric structures
built on it \cite{Woro89,Sita94-groups,DMH94-groups}. A finite
group together with a (bicovariant) differential calculus may be
regarded as a `finite Lie group'. The purpose of the present
paper is to develop differential geometry on such spaces.
In discrete (field) theories, discrete groups may appear as gauge
groups, as isometry groups, and as structures underlying
discrete space-time models. For example, the hypercubic lattice
underlying ordinary lattice (gauge) theories can be regarded as the
abelian group $\Ir^n$ (respectively, $\Ir^n_N$ with a positive integer
$N$, for a finite lattice). Lattice gauge theory can be understood as
gauge theory on this group with a bicovariant differential calculus
\cite{DMHS93}. Another example which fits into our framework is the
two-point-space used in \cite{Conn+Lott90} to geometrize models of
particle physics (see also \cite{CFF93,Sita94-groups,DMH94-groups}).
In this model the group $\Ir_2$ appears with the universal differential
calculus (which is bicovariant).
\vskip.2cm

Section \ref{Diffcalc} introduces to differential calculus on finite
groups and recalls the notions of left-, right- and bicovariance (using
the language of Hopf algebras).
For each bicovariant first order differential calculus on a finite
group and, more generally on a Hopf algebra, there is an operator
which acts on the tensor product of 1-forms and satisfies the
braid relation \cite{Woro89}. For a commutative finite group this is
simply the permutation operator, but less trivial structures arise in
the case of noncommutative groups. The generalized permutation operator
can be used to define symmetric and antisymmetric tensor fields. All
this is the subject of section \ref{Bm_iso}.
Linear connections on finite groups and corresponding invariance
conditions are considered in section \ref{Lincon}. Of particular
interest are linear connections which can be extended to tensor products
of 1-forms. We explore the restrictions on linear connections which
arise from the extension property. Appendix A more generally addresses
the problem of extending connections on two $\cal A$-bimodules, with
$\cal A$ any associative algebra, to a connection on their tensor
product.
In section \ref{Vecfield} we introduce vector fields on
finite groups and briefly discuss a possible concept of a metric.
In order to formulate, for example, metric-compatibility of a linear
connection, the concept of the `dual' of a connection is needed and the
problem of its extensibility (in the sense mentioned above) has to be
clarified. This is done in Appendix B for an arbitrary associative
algebra $\cal A$ and a connection on an $\cal A$-bimodule.
An example of a noncommutative finite group is elaborated in
section \ref{NgS3}. Section \ref{Concl} contains further discussion and
conclusions. Appendix C recalls how coactions on two bimodules
extend to a coaction on their tensor product. It is then shown that
extensions of invariant connections are again invariant connections
and that invariance is also carried over to the dual of a connection.
In Appendix D we briefly explore the concept of a `two-sided connection'
on a bimodule.
Appendix E deals with invariant tensor fields on finite groups.
Although in this work we concentrate on the case of finite sets
supplied with a group structure, in appendix F we indicate how
the formalism can be extended to the more general case of a finite
set with a finite group acting on it.
\vskip.2cm

Though originated from the development of `differential geometry'
on finite groups, some of our results are more general, they apply
to arbitrary associative algebras, respectively, Hopf algebras. We
therefore decided to separate them from the main part of the paper
and placed them into a series of appendices (A-D).

\section{Differential calculi on finite groups}
\label{Diffcalc}
\setcounter{equation}{0}
Every finite set can be supplied with a group structure.
If the number $N$ of elements is prime, then the only irreducible
group structure is $\Ir_N$, the additive abelian group of integers
modulo $N$. Differential calculi on discrete groups have been studied in
\cite{Sita94-groups,DMH94-groups}. More generally, differential calculus
on discrete sets has been developed in \cite{DMH94-graphs,DMHV}.
\vskip.2cm

Let $\cal A$ be the set of $\Cx$-valued functions on a finite set $G$.
With each element $g \in G$ we associate a function $e_g \in
{\cal A}$ via $e_g(g') = \delta_{g,g'}$. Then $e_g \, e_{g'} =
\delta_{g,g'} \, e_g$ and $\sum_{g \in G} e_g = \idty$
where $\idty$ is the unit in ${\cal A}$. Every function $f$ on
$G$ can be written as $f = \sum_{g \in G} f_g \, e_g$ with $f_g \in
\Cx$. Choosing a group structure on $G$, the latter induces a
{\em coproduct} $\Delta \, : \, {\cal A} \rightarrow {\cal A} \otimes
{\cal A}$ via\footnote{Here we make use of the fact that a group
defines a Hopf algebra. This formalism is adequate, in particular,
if one has in mind to generalize the structures considered in the
present work to noncommutative Hopf algebras like quantum groups.}
\be
     \Delta(f) (g,g') = f(g g')  \, .
\ee
In particular,
\be
     \Delta(e_g) =  \sum_{h \in G} e_h \otimes e_{h^{-1} g}  \, .
\ee
A {\em differential calculus} on $\cal A$ is an extension of
${\cal A}$ to a differential algebra $(\Omega, \d)$. Here
$\Omega = \bigoplus_{r=0}^\infty \Omega^r$ is a graded associative
algebra where $\Omega^0 = \cal A$. $\Omega^{r+1}$ is generated as an
$\cal A$-bimodule via the action of a linear operator $\d \, : \,
\Omega^r \rightarrow \Omega^{r+1}$ satisfying $\d^2=0$, $\d \idty =0$,
and the graded Leibniz rule
$\d(\varphi \varphi') = (\d \varphi) \, \varphi'+(-1)^r \varphi \, \d
\varphi'$ where $\varphi \in \Omega^r$.\footnote{In \cite{DMHV} a
discrete set together with a differential calculus on it has been
called {\em discrete differential manifold}. This notion was motivated
by a far-reaching analogy \cite{DMH94-graphs} with the continuum case
where $\Omega$ is the algebra of differential forms on a manifold and
$\d$ the exterior derivative.}
\vskip.2cm

It is convenient \cite{DMH94-graphs,DMH94-groups} to introduce the
special 1-forms
\begin{eqnarray}
     e_{g,g'} := e_g \, \d e_{g'}    \qquad  (g \neq g') \, , \qquad
     e_{g,g} := 0
\end{eqnarray}
and the $(r-1)$-forms
\begin{eqnarray}            \label{e...}
  e_{g_1, \ldots, g_r} := e_{g_1, g_2} \, e_{g_2, g_3} \cdots
                        e_{g_{r-1}, g_r} \;.
\end{eqnarray}
They satisfy
\begin{eqnarray}       \label{e-e}
   e_{g_1, \ldots, g_r} \, e_{h_1, \ldots, h_s}
   = \delta_{g_r, h_1} \, e_{g_1, \ldots, g_r, h_2, \ldots, h_s} \; .
\end{eqnarray}
The operator $\d$ acts on them as follows,
\begin{eqnarray}
  \d e_{g_1, \ldots, g_r}
    =  \sum_{h \in G} \lbrack e_{h, g_1, \ldots, g_r}
      - e_{g_1, h, g_2, \ldots, g_r}
      + e_{g_1, g_2, h, g_3, \ldots, g_r}
      - \ldots + (-1)^r \, e_{g_1, \ldots, g_r, h} \rbrack \; .
                                       \label{de...}
\end{eqnarray}
If no further relations are imposed, one is dealing with
the `universal differential calculus' $(\tilde{\Omega}, \tilde{\d})$.
The $e_{g_1, \ldots, g_r}$ with $g_i \neq g_{i+1}$ ($i=1, \ldots,
r-1$) then constitute a basis over $\Cx$ of
$\tilde{\Omega}^{r-1}$ for $r>1$ \cite{DMH94-graphs}.
Every other differential calculus on $G$ is obtained from
$\tilde{\Omega}$ as the quotient with respect to some two-sided
differential ideal. Up to first order, i.e., the level of 1-forms,
every differential calculus on $G$ is obtained by setting some of the
$e_{g,g'}$ to zero. Via (\ref{e...}) and (\ref{de...}) this induces
relations for forms of higher grade. In addition, or alternatively, one
may also factor out ideals generated by forms of higher grade.
Every first order differential calculus on $G$ can be described
by a (di)graph the vertices of which are the elements of $G$
and there is an arrow pointing from a vertex $g$ to a vertex $g'$ iff
$e_{g,g'} \neq 0$ (see \cite{DMH94-graphs} for further details).
\vskip.2cm

A differential calculus on $G$ (or, more generally, any Hopf algebra
$\cal A$) is called {\em left-covariant}
\cite{Woro89} if there is a linear map $\ldel \, : \, \Omega^1
\rightarrow {\cal A} \otimes \Omega^1$ such that
\be
  \ldel (f \, \varphi \, f') = \Delta(f) \, \ldel (\varphi) \,
  \Delta(f') \qquad  \forall f,f' \in {\cal A}, \, \varphi \in \Omega^1
\ee
and
\be
      \ldel \circ \d = (id \otimes \d) \circ \Delta \; .
\ee
As a consequence,\footnote{Without refering to the Hopf-algebraic
language, left-covariance of a differential calculus basically means
${\cal L}_g \d = \d {\cal L}_g$ where ${\cal L}_g$ is the action on
$\cal A$ induced by the left multiplication of elements of $G$ by
$g$, see (\ref{left-action}). Then ${\cal L}_g e_{h,h'} = e_{g^{-1} h ,
g^{-1} h'}$ which corresponds to (\ref{e-lcoaction}).}
\be                       \label{e-lcoaction}
   \ldel (e_{g,g'}) = \sum_{h \in G} e_{h^{-1}} \otimes e_{h g,
                      h g'} \; .
\ee
Hence, in order to find the left-covariant differential calculi
on $G$, we have to determine the orbits of all elements of
$(G \times G)'$ where the prime indicates omission of the diagonal
(i.e., $(G \times G)' = (G \times G) \setminus \lbrace (g,g) \mid g \in
G \rbrace$) with respect to the left action $(g,g') \mapsto
(h g, h g')$. In the graph picture, left-covariant first order
differential calculi are obtained from the universal one (which is
left-covariant) by deleting corresponding orbits of arrows.
\vskip.2cm

For a left-covariant\footnote{All following formulas and statements
which do not make explicit reference to a coaction are actually
valid without the assumption of left-covariance. The second formula
in (\ref{lMC}) is based on it, however, and it is this invariance
property which justifies our definition in (\ref{lMC}).}
differential calculus, there are left-invariant
`Maurer-Cartan' 1-forms \cite{Sita94-groups,DMH94-groups}
\be
   \theta^g := \sum_{h \in G} e_{hg,h}
             = e_{\un{h} \, g,\un{h}} = e_{\un{h},\un{h} \, g^{-1}}
   \, , \quad
   \ldel (\theta^g) = \idty \otimes \theta^g  \; .
                        \label{lMC}
\ee
Here we have introduced a summation convention. If an index is
underlined, this means summation over all group elements. Note
that $\theta^e = 0$ according to the above definition. Furthermore,
the Maurer-Cartan forms with $g \neq e$ are in
one-to-one correspondence with left orbits in $(G \times G)'$.
All left-covariant differential calculi (besides the universal one) are
therefore obtained by setting some of the $\theta^g$ (of the
universal calculus) to zero. The nonvanishing $\theta^g$ then
constitute a left (or right) $\cal A$-module basis for $\Omega^1$
since
\be                                      \label{e-theta}
    e_{h,g} = e_h \, \theta^{g^{-1}h} = \theta^{g^{-1}h} \, e_g \; .
\ee
As a generalization of the last equality we have the simple commutation
relations
\be                            \label{f-theta}
         f \, \theta^g = \theta^g \, {\cal R}_g f
\ee
where ${\cal R}_g$ denotes the action of $G$ on $\cal A$ induced by
right multiplication, i.e.,
\be
  ({\cal R}_g f)(h) := f(hg) \quad  (\forall f \in {\cal A})
  \; , \quad
   {\cal R}_g \, {\cal R}_h = {\cal R}_{gh}  \; .
                            \label{right-action}
\ee
\vskip.2cm

The equation (\ref{e-theta}) can be used to prove the Maurer-Cartan
equations
\be
   \d \theta^h = - C^h{}_{\un{g},\un{g}'} \, \theta^{\un{g}'}
                 \, \theta^\un{g}
\ee
with the `structure constants'
\be                         \label{struct_const}
  C^h{}_{g,g'} := - \delta^h_g - \delta^h_{g'} + \delta^h_{g g'} \; .
\ee
These have the property
\be
   C^{ad(h)g}{}_{ad(h)g',ad(h)g''} = C^g{}_{g',g''}
   \qquad   \forall  h \in G              \label{C-adinv}
\ee
where $ad$ denotes the adjoint action of $G$ on $G$, i.e., $ad(h)g =
hgh^{-1}$.
\vskip.2cm

A differential calculus on $G$ is called {\em right-covariant} if there
is a linear map $\rdel \, : \, \Omega^1 \rightarrow \Omega^1 \otimes
{\cal A}$ such that
\be
 \rdel (f \, \varphi \, f') = \Delta(f) \, \rdel(\varphi) \, \Delta(f')
 \, , \quad  \rdel \circ \d = (\d \otimes id) \circ \Delta \; .
\ee
This implies
\be
 \rdel (e_{g,g'}) = \sum_{h \in G} e_{g h, g' h} \otimes e_{h^{-1}} \; .
\ee
For a right-covariant differential calculus, there are right-invariant
Maurer-Cartan 1-forms,
\be
   \omega^g := e_{g \, \un{h},\un{h}}
               \; , \qquad
   \rdel (\omega^g) = \omega^g \otimes \idty  \; ,
\ee
which satisfy
\be
      e_{h,g} &=& e_h \, \omega^{hg^{-1}} = \omega^{hg^{-1}} \, e_g \\
  \d \omega^h &=& - C^h{}_{\un{g},\un{g}'} \, \omega^\un{g}
                   \, \omega^{\un{g}'}
\ee
and
\be
           f \, \omega^g = \omega^g \, {\cal L}_g f
\ee
where ${\cal L}_g$ denotes the action of $G$ on $\cal A$ induced by
left multiplication,
\be
   ({\cal L}_g f)(h) := f(gh) \quad  (\forall f \in {\cal A})
   \; , \quad
   {\cal L}_g \, {\cal L}_h = {\cal L}_{hg} \; .
                                            \label{left-action}
\ee
\vskip.2cm

A differential calculus is {\em bicovariant} if it is left- and
right-covariant. Then, in the case under consideration,
\be
    \rdel (\theta^g) &=& \theta^{ad(\un{h})g} \otimes e_{\un{h}} \; .
\ee
It follows that bicovariant calculi are in one-to-one
correspondence with unions of conjugacy classes different from
$\lbrace e \rbrace$. Obviously,
\be
    \rho := \theta^\un{g} = \omega^\un{g} = e_{\un{g},\un{g}'}
\ee
is a {\em bi-invariant} 1-form.
\vskip.2cm

In the following we list some useful formulas. For $f \in {\cal A}$
we find
\be                           \label{df-rho}
   \d f = \lbrack \rho , f \rbrack
        = (\ell_\un{g} f) \; \theta^\un{g}
        = (r_\un{g} f) \; \omega^\un{g}
\ee
where
\be
    \ell_g f := {\cal R}_{g^{-1}} f - f \; ,  \qquad
       r_g f := {\cal L}_{g^{-1}} f - f   \; .
\ee
Using (\ref{left-action}) and (\ref{right-action}), it is easy to check
that
\be
   \ell_g \, \ell_{g'} = C^\un{h}{}_{g',g} \; \ell_\un{h} \; ,  \qquad
   r_g \, r_{g'} = C^\un{h}{}_{g,g'} \; r_\un{h} \; .
\ee
The 1-forms $\theta^g$ and $\omega^g$ are related as follows,
\be
   \theta^g = e_{\un{h}} \, \omega^{ad(\un{h}) g} \; , \quad
   \omega^g = e_\un{h} \, \theta^{ad(\un{h}^{-1}) g}     \; .
\ee
\vskip.2cm

In the following sections we restrict our considerations to
differential calculi which are at least left-covariant. As
already mentioned, in this case the set of nonvanishing left-invariant
Maurer-Cartan 1-forms $\theta^g$ is a basis of $\Omega^1$ as a left
$\cal A$-module. It is then convenient to introduce the subset
$\hat{G} := \lbrace g \in G \mid \theta^g \neq 0 \rbrace$ of $G$.
If not said otherwise, indices will be restricted to $\hat{G}$ in what
follows. This does {\em not} apply to our summation convention,
however. Underlining an index still means summation over all elements
of $G$, though in most cases the sum reduces to a sum over $\hat{G}$
(but see (\ref{sig-theta}) for an exception).

\section{The canonical bimodule isomorphism for a bicovariant
         differential calculus}
\label{Bm_iso}
\setcounter{equation}{0}
For a bicovariant differential calculus there is a unique bimodule
isomorphism $\sigma \, : \, \Omega^1 \otimes_{\cal A} \Omega^1
\rightarrow \Omega^1 \otimes_{\cal A} \Omega^1$ such that
\be                \label{sig_th_om}
   \sigma (\theta \otimes_{\cal A} \omega) = \omega \otimes_{\cal A}
                                             \theta
\ee
for all left-invariant 1-forms $\theta$ and right-invariant 1-forms
$\omega$ \cite{Woro89}.
We have $\sigma (\rho \otimes_{\cal A} \rho) = \rho
\otimes_{\cal A} \rho$ since $\rho$ is bi-invariant.
Furthermore,
\be
     \sigma (\theta^g \otimes_{\cal A} \theta^{g'})
 &=& \sigma (\theta^g \otimes_{\cal A} e_{\un{h}} \,
     \omega^{ad(\un{h}) g'})                               \nonumber \\
 &=& \sigma (\theta^g \, e_{\un{h}} \otimes_{\cal A}
     \omega^{ad(\un{h}) g'})                               \nonumber \\
 &=& e_{\un{h} \, g} \, \sigma (\theta^g \otimes_{\cal A}
     \omega^{ad(\un{h}) g'})                               \nonumber \\
 &=& e_{\un{h} \, g} \, \omega^{ad(\un{h}) g'} \otimes_{\cal A} \theta^g
                                                           \nonumber \\
 &=& e_{\un{h} \, g} e_{\un{h}'} \, \theta^{ad({\un{h}'}^{-1}
     \un{h}) g'} \otimes_{\cal A} \theta^g                 \nonumber \\
 &=& e_\un{h} \, \theta^{ad(g^{-1}) g'} \otimes_{\cal A} \theta^g
                                           \label{sig-theta}
\ee
which implies
\be
    \sigma (\theta^g \otimes_{\cal A} \theta^{g'}) = \theta^{ad(g^{-1})
    g'} \otimes_{\cal A} \theta^g   \, .    \label{sigma-theta}
\ee
In particular, $\sigma (\theta^g \otimes_{\cal A} \theta^g) = \theta^g
\otimes_{\cal A} \theta^g$.
More generally, it is possible to calculate an expression for higher
powers of $\sigma$. By induction one can prove that
\be
    \sigma^{2n-1} (\theta^g \otimes_{\cal A} \theta^{h})
   &=& \theta^{ad(g^{-1} h^{-1})^{n}h} \otimes_{\cal A}
       \theta^{ad(g^{-1}h^{-1})^{n-1}g}            \label{sig2n-1} \\
    \sigma^{2n} (\theta^g \otimes_{\cal A} \theta^{h})
   &=& \theta^{ad(g^{-1} h^{-1})^{n}g} \otimes_{\cal A}
       \theta^{ad(g^{-1}h^{-1})^nh}                \label{sig2n}
\ee
for all $n \geq 1$. With the help of the last formula one arrives at
the following result.
\vskip.3cm
\noindent
{\em Proposition 3.1.} For a finite group $G$ and a bicovariant
first order differential calculus on it, the associated bimodule
isomorphism $\sigma$ satisfies
\be                     \label{sigid}
           \sigma^{2 |ad(G)|} = id
\ee
where $|ad(G)|$ denotes the number of elements of $ad(G) := \lbrace
ad(g) \mid g \in G \rbrace$, the group of inner automorphisms of $G$.
\vskip.1cm \noindent
{\em Proof:} For $a \in ad(G)$ let $\langle a \rangle$ denote the
cyclic subgroup of $ad(G)$ generated by $a$. Since
$ad(G)$ is a finite group, $|\langle a \rangle|$ is finite and
$a^{| \langle a \rangle |} = id$. Furthermore, $| \langle a \rangle |$
is a divisor of $|ad(G)|$ by Lagrange's theorem. Now (\ref{sigid})
follows from (\ref{sig2n}).                  \hfill  {\Large $\Box$}
\vskip.3cm
\noindent
We define the {\em order} of $\sigma$ as the smallest positive
integer $m$ such that $\sigma^m = id$ and denote it as $|\sigma|$.
The previous proposition then tells us that $|\sigma| \leq 2 \,
|ad(G)|$. Our next result shows that, in general, equality does
not hold. For the symmetric groups ${\cal S}_n$ with $n>3$ one
finds that $|\sigma| < 2 \, |ad({\cal S}_n)|$.
\vskip.3cm
\noindent
{\em Proposition 3.2.} For the symmetric group ${\cal S}_n$, $n \geq
3$, with the universal first order differential calculus, we
have\footnote{Another way to describe the number on the rhs is the
following. Write all the factors $2,\ldots,n$ of $|{\cal S}_n| = n!$
as products of powers of primes. Then $|\sigma|$ is twice
the product of all different primes, each taken to the power which
is the highest with which the prime appears in the set of factors
$2,\ldots,n$.}
\be
   |\sigma| = 2 \, n \, \prod_{k=1}^{n-2} {n-k \over
   \mbox{gcd} \lbrack n (n-1) \cdots (n-k+1), n-k \rbrack }
\ee
where $\mbox{gcd} \lbrack \ell, \ell' \rbrack$ denotes the greatest
common divisor of positive integers $\ell$ and $\ell'$.
\vskip.1cm
\noindent
{\em Proof:}
For $n > 2$ the center of ${\cal S}_n$ is trivial \cite{Pass68} and
the group of inner automorphisms is therefore isomorphic with
${\cal S}_n$ itself. Every element $g \in {\cal S}_n$, $g \neq e$, can
be written as a product of disjoint (and thus commuting) cycles. Hence
$g^\ell = e$ with $\ell := n \, \prod_{k=1}^{n-2} (n-k)/ \mbox{gcd}
\lbrack n (n-1) \cdots (n-k+1), n-k \rbrack$. By construction,
$\ell$ is the smallest positive integer with the property $g^{\ell} = e$
for all $g \in {\cal S}_n$, since for each divisor of $\ell$ there is a
cyclic subgroup of order equal to this divisor in ${\cal S}_n$ (given by
a cycle of length equal to the divisor).
For the universal differential calculus on ${\cal S}_n$, the
statement in the proposition now follows from (\ref{sig2n-1}) and
(\ref{sig2n}).
\hspace*{1cm}                               \hfill  {\Large $\Box$}
\vskip.2cm

For a bicovariant differential calculus, it is natural to consider the
following symmetrization and antisymmetrization operators acting on
$\Omega^1 \otimes_{\cal A} \Omega^1$,
\be
   {\bf S} := {1 \over 2} (id + \sigma) \, , \quad
   {\bf A} := {1 \over 2} (id - \sigma) \; .
\ee
In general, $\sigma^2 \neq id$, so that these are not projections. It is
therefore not quite straightforward how to define symmetry and
antisymmetry for an element $\alpha \in \Omega^1 \otimes_{\cal A}
\Omega^1$. We suggest the following notions,
\begin{eqnarray*}
  \alpha \; \mbox{is w-symmetric iff}& & \alpha \in
                {\bf S}(\Omega^1 \otimes_{\cal A} \Omega^1)
                = \mbox{im} \, {\bf S} \\
  \alpha \; \mbox{is s-symmetric iff}& & {\bf S}(\alpha) = \alpha  \\
  \alpha \; \mbox{is w-antisymmetric iff}& & \alpha \in
                {\bf A}(\Omega^1 \otimes_{\cal A} \Omega^1)
                = \mbox{im} \, {\bf A} \\
  \alpha \; \mbox{is s-antisymmetric iff}& & {\bf A}(\alpha) = \alpha
\end{eqnarray*}
where `w' and `s' stand for `weakly' and `strongly',
respectively.\footnote{The conditions ${\bf S}(\alpha) = \alpha$ and
${\bf A}(\alpha) = \alpha$ are equivalent to $\alpha\in \ker {\bf A}$
and $\alpha\in \ker {\bf S}$, respectively.}
Examples are treated in section \ref{NgS3}.1 and Appendix E.
The notions of s-symmetry and w-antisymmetry are complementary
in the following sense.
\vskip.3cm
\noindent
{\em Proposition 3.3.} For each bicovariant differential calculus
$\Omega$ on a finite group $G$, the space $\Omega^1 \otimes_{\cal A}
\Omega^1$ decomposes into direct sums
\be                                \label{dirsums}
    \Omega^1\otimes_{\cal A}\Omega^1
 = \ker {\bf A} \oplus \mbox{im} \, {\bf A}
 = \ker {\bf S} \oplus \mbox{im} \, {\bf S} \; .
\ee
{\em Proof:} In order to show the first equality, it is sufficient to
prove that $\ker {\bf A} \cap \mbox{im} \, {\bf A}
=0$ since $\Omega^1 \otimes_{\cal A} \Omega^1$ is a
finite-dimensional vector space over $\Cx$ and $\bf A$ a linear map.
Let $\alpha$ be an element of $\ker {\bf A} \cap \mbox{im} \, {\bf A}$.
Then $\sigma(\alpha) = \alpha$ and $\alpha
= (\sigma - id)(\beta)$ with an element $\beta \in {\bf A}(\Omega^1
\otimes_{\cal A} \Omega^1)$. Using (\ref{sigid}) we obtain
\begin{eqnarray*}
     0 = (\sigma^{2 \mid ad(G) \mid}-id)(\beta)
       = (\sum_{k=0}^{2 \mid ad(G) \mid-1}\sigma^k)(\sigma-id)(\beta)
       = (\sum_{k=0}^{2 \mid ad(G) \mid - 1} \sigma^k)(\alpha)
       = 2 \, |ad(G)| \, \alpha \; .
\end{eqnarray*}
Hence, $\alpha = 0$. In the same way, the second equality in
(\ref{dirsums}) is varified with the help of
$0=(\sum_{k=0}^{2 \mid ad(G) \mid - 1}(-1)^{k+1} \sigma^k)(\sigma+id)$.
                                            \hfill  {\Large $\Box$}

\vskip.3cm

We note that $\sigma$ satisfies the braid equation\footnote{The fact
that $\sigma$ satisfies the braid relation has the following origin.
Let $V$ be a vector space and $\Phi$ a map $V \rightarrow \mbox{End}
(V)$.  The map $\tilde{\sigma} \, : \, V \otimes V \rightarrow
V \otimes V$ defined by $\tilde{\sigma}(x \otimes y) := y \otimes
\Phi_y x$ then satisfies the braid equation if and only if
$\Phi_x \circ \Phi_y = \Phi_{\Phi_x y} \circ \Phi_x$ for all
$x,y \in V$. In particular, if $V$ is the group algebra of a (not
necessarily finite) group $G$, then $\Phi = ad$ satisfies this
equation.}
\be
     (id \otimes \sigma) (\sigma \otimes id) (id \otimes \sigma)
   = (\sigma \otimes id) (id \otimes \sigma) (\sigma \otimes id)
\ee
(see also \cite{Woro89}).
In \cite{Woro89} Woronowicz implemented a generalized wedge product
by taking the quotient of $\Omega^1 \otimes_{\cal A} \Omega^1$
with respect to the subbimodule of s-symmetric tensors, i.e.,
$\Omega^2 = (\Omega^1 \otimes_{\cal A} \Omega^1)/ \ker {\bf A}$
which can be identified with the space of w-antisymmetric
tensors.\footnote{Alternatively, one may think of
implementing a generalized wedge product by taking the quotient
with respect to w-symmetric tensors. However, this turns out to be too
restrictive, in general (see section 6.1 and the example in Appendix E).
Moreover, one may also consider corresponding (anti)symmetry conditions
obtained from those given above by replacing $\sigma$ by some power of
$\sigma$ and use them to define a wedge product. These possibilities
reflect the fact that there are several differential algebras with the
same space of 1-forms. These define different discrete differential
manifolds \cite{DMH94-graphs,DMHV}. The choice made by Woronowicz is
uniquely distinguished by the property that bicovariance extends to the
whole differential algebra, see Theorem 4.1 in \cite{Woro89}.}
Then $\d \rho = \rho^2 = 0$.
\vskip.2cm
\noindent
{\em Example.} We consider a set of three elements with the group
structure $\Ir_3$. The $\Ir_3$ left-covariant first order differential
calculi on this set are then represented by the graphs in Fig. 1.

\unitlength1.cm
\begin{picture}(12.,1.8)(-2.8,-0.3)
\thicklines
\put(0.,1.) {\circle*{0.1}}
\put(1.,0.) {\circle*{0.1}}
\put(-1.,0.) {\circle*{0.1}}
\put(-0.8,0.05) {\vector(1,0){1.6}}
\put(0.8,-0.07) {\vector(-1,0){1.6}}
\put(0.8,0.1) {\vector(-1,1){0.7}}
\put(0.2,0.9) {\vector(1,-1){0.7}}
\put(-0.1,0.8) {\vector(-1,-1){0.7}}
\put(-0.9,0.2) {\vector(1,1){0.7}}
\put(3.,1.) {\circle*{0.1}}
\put(4.,0.) {\circle*{0.1}}
\put(2.,0.) {\circle*{0.1}}
\put(3.9,0.) {\vector(-1,0){1.8}}
\put(3.1,0.9) {\vector(1,-1){0.8}}
\put(2.1,0.1) {\vector(1,1){0.8}}
\put(6.,1.) {\circle*{0.1}}
\put(7.,0.) {\circle*{0.1}}
\put(5.,0.) {\circle*{0.1}}
\put(5.1,0.) {\vector(1,0){1.8}}
\put(6.9,0.1) {\vector(-1,1){0.8}}
\put(5.9,0.9) {\vector(-1,-1){0.8}}
\put(9.,1.) {\circle*{0.1}}
\put(10.,0.) {\circle*{0.1}}
\put(8.,0.) {\circle*{0.1}}
\end{picture}

\begin{center}
\begin{minipage}{10.5cm}
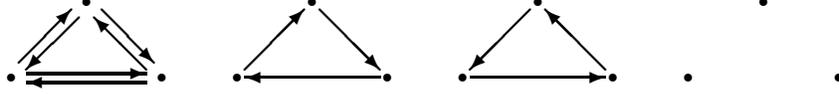
\centerline{\bf Fig. 1}
\vskip.1cm
\noindent
The digraphs which determine all left-covariant first order
differential calculi on $\Ir_3$.
\end{minipage}
\end{center}
\vskip.2cm
\noindent
The last of these graphs has no arrows and corresponds to the trivial
differential calculus (where $\d \equiv 0$). For a commutative group,
as in the present example, bicovariance does not lead to additional
conditions. Since left- and right-invariant Maurer-Cartan forms
coincide for a commutative group, the map $\sigma$ acts on
left-invariant forms simply as permutation, i.e., $\sigma
(\theta^g \otimes_{\cal A} \theta^{g'}) = \theta^{g'} \otimes_{\cal A}
\theta^g$. In particular, $\sigma^2 = id$ in accordance with
Proposition 3.1 and the wedge product determined by $\sigma$ is
therefore the ordinary one for left-invariant 1-forms, though we still
do not have anticommutativity of the product of two 1-forms, in general.
                                      \hfill  {\large \bf $\Box$}
\vskip.2cm
\noindent
More complicated maps $\sigma$ with $\sigma^2 \neq id$ arise from a
noncommutative group structure. In section \ref{NgS3} we elaborate in
some detail the case of the symmetric group ${\cal S}_3$.

\section{Linear connections on a finite group}
\label{Lincon}
\setcounter{equation}{0}
Let $\cal A$ be an associative algebra and $\Omega$ a
differential calculus on it. A {\em connection} on a left
$\cal A$-module $\Gamma$ is a map\footnote{Similarly, a connection
on a right $\cal A$-module is a map $\nabla \, : \, \Gamma
\rightarrow \Gamma \otimes_{\cal A} \Omega^1$ with $\nabla
(\gamma \, f) = (\nabla \gamma) \, f + \gamma \otimes_{\cal A} \d f$.
A left (right) module over an associative algebra $\cal A$ has a
connection with respect to the universal first order differential
calculus if and only if it is projective \cite{Cunt+Quil95} (see also
\cite{Conn85}).}
\begin{eqnarray}
 \nabla \; : \; \Gamma \rightarrow \Omega^1 \otimes_{\cal A}
           \Gamma                     \label{connection}
\end{eqnarray}
such that
\begin{eqnarray}
  \nabla (f \, \gamma) = \d f \otimes_{\cal A} \gamma
  + f \, \nabla \gamma
  \qquad \quad \forall f \in {\cal A}, \, \gamma \in \Gamma
                           \label{conn-prop}
\end{eqnarray}
\cite{Conn85}. A connection on $\Gamma$ can be extended to a map
$\Omega \otimes_{\cal A} \Gamma \rightarrow \Omega \otimes_{\cal A}
\Gamma$ via
\be
  \nabla(\varphi \Psi) = (\d \varphi) \Psi + (-1)^r \varphi \nabla \Psi
\ee
for $\varphi \in \Omega^r$ and $\Psi \in \Omega \otimes_{\cal A} \Gamma$.
Then $\nabla^2$, which is a left $\cal A$-module homomorphism, defines
the {\em curvature} of the connection.
\vskip.2cm

In the following we consider the particular case where $\Gamma =
\Omega^1$, the space of 1-forms of a differential calculus on $\cal A$.
A connection is then called a {\em linear connection}. It is a map
$\nabla \, : \, \Omega^1 \rightarrow \Omega^1 \otimes_{\cal A}
\Omega^1$.\footnote{This pretends that there should be a kind of
symmetry with respect to the two factors of the tensor product. From
the general formula (\ref{connection}) it should be clear that the
two factors play very different roles.}
The {\em torsion} of a linear connection may be defined as
\be
       T = \d - \pi \circ \nabla
\ee
where $\pi$ is the projection $\Omega^1 \otimes_{\cal A}
\Omega^1 \rightarrow \Omega^2 \,$.\footnote{For a bicovariant
first order differential calculus on a Hopf algebra, the
choice made in \cite{Woro89} is $\pi = {\bf A}$, see section
\ref{Bm_iso}.}
\vskip.2cm

If $\Omega^1$ has a left $\cal A$-module basis $\theta^i, \,
i=1,\ldots,n$, the action of a linear connection on a 1-form $\varphi
= \varphi_i \, \theta^i$ (summation convention) is given by
\begin{eqnarray}
    \nabla \varphi = D \varphi_i \otimes_{\cal A} \theta^i
\end{eqnarray}
where we have introduced
\begin{eqnarray}
    D \varphi_i := \d \varphi_i - \varphi_j \, \omega^j{}_i
\end{eqnarray}
with connection 1-forms
\begin{eqnarray}                       \label{lin_conn}
   \omega^i{}_j = \Gamma^i{}_{j k} \; \theta^k
\end{eqnarray}
defined by
\begin{eqnarray}                  \label{conn-1forms}
   \nabla \theta^i = - \omega^i{}_j \otimes_{\cal A} \theta^j   \; .
\end{eqnarray}
Extending a linear connection on $\Omega^1$ to a connection on
$\Omega \otimes_{\cal A} \Omega^1$,
\be
 \nabla^2 \theta^i = - \Omega^i{}_j \otimes_{\cal A} \theta^j
\ee
defines curvature 2-forms\footnote{More precisely, these are the
components of the curvature with respect to the (arbitrary) left
$\cal A$-module basis $\theta^i$.}
for which we obtain the familiar formula
\be
  \Omega^i{}_j = \d \omega^i{}_j + \omega^i{}_k \,
  \omega^k{}_j    \; .
\ee
\vskip.3cm
\noindent
{\em Remark.} Under a change of basis $\theta^i \mapsto a^i{}_j \,
\theta^j$ where $a$ is an invertible matrix with entries in $\cal A$,
we have the tensorial transformation properties $\varphi_i \mapsto
\varphi_j \, (a^{-1})^j{}_i$ and $D\varphi_i \mapsto D \varphi_j \,
(a^{-1})^j{}_i$. For the connection 1-forms and the curvature 2-forms
one finds the familiar transformation laws $\omega^i{}_j \mapsto a^i{}_k
\, \omega^k{}_l \, (a^{-1})^l{}_j + a^i{}_k \, \d (a^{-1})^k{}_j \,$
and $\, \Omega^i{}_j \mapsto a^i{}_k \, \Omega^k{}_l \, (a^{-1})^l{}_j$.
It should be noticed, however, that the components $R^i{}_{jkl}$ of
the 2-forms $\Omega^i{}_j$ with respect to the generators $\theta^k
\theta^l$ of $\Omega^2$ do {\em not} transform in this simple way
if functions (here the entries of the transformation matrix $a$) do not
commute with all 1-forms (here the basis 1-forms $\theta^k$), as in the
case of a differential calculus on a finite set.\footnote{On the other
hand, it is precisely such a change of the ordinary transformation law
for components of forms which underlies the derivation of the lattice
gauge theory action in \cite{DMHS93}.}
                                             \hfill {\Large $\Box$}
\vskip.3cm

Now we turn to the special case of a left-covariant differential
calculus on a finite group. As a left $\cal A$-module basis we choose
the set of left-invariant Maurer-Cartan 1-forms $\theta^g$. Except where
stated otherwise, indices are restricted to $\hat{G}$. For $\varphi =
\varphi_\un{g} \, \theta^\un{g}$ one finds
\be
  \nabla \varphi = ({\cal R}_{\un{g}^{-1}} \varphi_{\un{g}'}
                   - \varphi_\un{h} \, U^\un{h}{}_{\un{g}',\un{g}}) \,
                   \theta^\un{g} \otimes_{\cal A} \theta^{\un{g}'}
\ee
where
\be
      U^h{}_{g',g} := \delta^h_{g'} + \Gamma^h{}_{g',g}  \; .
\ee
As a consequence,
\be
  \nabla \varphi = 0       \quad \Leftrightarrow \quad
    {\cal R}_{g^{-1}} \varphi_{g'}
  = \varphi_\un{h} \, U^\un{h}{}_{g',g} \;.
\ee
Left-invariance of $\nabla$ (see Appendix C) is equivalent to
$\Gamma^g_{g',g''} \in \Cx$. If the differential calculus is
bicovariant, evaluation of the right-invariance condition (see
Appendix C) then leads to\footnote{Note that, for a bicovariant
differential calculus, $ad(h)g \in \hat{G}$ whenever $g \in \hat{G},
h \in G$.}
\be
   \Gamma^{ad(h)g}_{ad(h)g', ad(h)g''} = \Gamma^g_{g',g''}
   \qquad \forall g,g',g'' \in \hat{G}, \; \forall h \in G  \; .
\ee
For a left-invariant connection,\footnote{For the universal differential
calculus this implies $\Omega = 0 \; \Leftrightarrow \; U_g \, U_{g'}
= U_{g g}$ where $U_g := (U^h{}_{h',g})$. The curvature thus measures
the deviation of the matrices $U_g$ from being a representation of the
group $G$.}
\be
   \Omega^g{}_{g'} = (\Gamma^g{}_{\un{g}'',\un{h}} \,
   \Gamma^{\un{g}''}{}_{g',\un{h}'} -
  C^{\un{g}''}{}_{\un{h}',\un{h}} \, \Gamma^g{}_{g',\un{g}''}) \,
  \theta^\un{h} \, \theta^{\un{h}'} \; .
\ee
Applying $T$ to the left-invariant basis, we find
\be
 T(\theta^h) = ( \Gamma^h{}_{\un{g}',\un{g}} - C^h{}_{\un{g}',\un{g}} )
               \, \theta^\un{g} \, \theta^{\un{g}'}   \; .
\ee
The constants $C^h{}_{g',g}$ are those defined in (\ref{struct_const}).
\vskip.3cm
\noindent
{\em Example 1.} For a bicovariant (first order) differential calculus,
\be                          \label{canon-conn}
   \nabla^\sigma \varphi :=  \rho \otimes_{\cal A} \varphi
   - \sigma(\varphi \otimes_{\cal A} \rho)
\ee
defines a linear left $\cal A$-module connection. Here $\sigma$ is the
canonical bimodule isomorphism. For this connection the Maurer-Cartan
1-forms $\theta^g$ are covariantly constant, i.e., $\nabla^\sigma
\theta^g = 0$.
As a consequence, the connection is bi-invariant, the curvature vanishes
and the torsion is given by $T(\theta^h) = - C^h{}_{\un{g}', \un{g}} \,
\theta^\un{g} \, \theta^{\un{g}'} = \d \theta^h$.       \\
The connection (\ref{canon-conn}) can be generalized to a family of
bi-invariant\footnote{The 1-form $\rho$ and the bimodule isomorphism
$\sigma$ are bi-invariant \cite{Woro89}. The bi-invariance of the
connections (\ref{conn-fam}) then follows from Proposition
\ref{Coact}.3.}
left $\cal A$-module connections,
\be                              \label{conn-fam}
   \nabla^{(\lambda_0,\ldots,\lambda_{|\sigma|-1})} \varphi :=
   \rho \otimes_{\cal A} \varphi - \sum_{n=0}^{|\sigma|-1} \lambda_n \,
   \sigma^n (\varphi \otimes_{\cal A} \rho)
\ee
where $|\sigma|$ is the order of $\sigma$ (see section \ref{Bm_iso})
and $\lambda_n \in \Cx$.\footnote{More generally, if $\Psi$
is any left $\cal A$-module homomorphism $\Omega^1 \otimes_{\cal A}
\Omega^1 \rightarrow \Omega^1 \otimes_{\cal A} \Omega^1$, then
$\nabla^\Psi \varphi := \rho \otimes_{\cal A} \varphi - \Psi(\varphi
\otimes_{\cal A} \rho)$ is a linear connection, see also
\cite{GMMM95}.}
It includes $\nabla^{\sigma^{-1}} := \nabla^{(0,\ldots,0,1)}$. With
respect to this connection the right-invariant Maurer-Cartan 1-forms
$\omega^g$ are covariantly constant, i.e., $\nabla^{\sigma^{-1}}
\omega^g = 0$, and the curvature also vanishes. The two connections
$\nabla^\sigma$ and $\nabla^{\sigma^{-1}}$ provide us with analogues of
the $(+)$- and $(-)$-parallelism on Lie groups (see \cite{Eise61},
\S 50). Corresponding {\em right} $\cal A$-module connections with
these properties are given by $\sigma^{-1} \circ \nabla^\sigma$ and
$\sigma \circ \nabla^{\sigma^{-1}}$.
                                     \hfill {\Large $\Box$}
\vskip.3cm
\noindent
{\em Example 2.} For a bicovariant differential calculus with the
(generalized) wedge product as defined by Woronowicz (see section
\ref{Bm_iso}), the torsion and the curvature of a left-invariant
linear (left $\cal A$-module) connection are given by
\be
   T(\theta^h) &=& {1 \over 2} \, ( \Gamma^h{}_{\un{g},\un{g}'}
   - \Gamma^h{}_{ad(\un{g}) \un{g}', \un{g}}
   - C^h{}_{\un{g},\un{g}'} + C^h{}_{ad(\un{g}) \un{g}', \un{g}} ) \,
   \theta^{\un{g}} \otimes_{\cal A} \theta^{\un{g}'}   \\
   \Omega^g{}_{g'} &=& {1 \over 2} \, \lbrack
   \Gamma^g{}_{\un{g}'',\un{h}} \, \Gamma^{\un{g}''}{}_{g',\un{h}'}
   - \Gamma^g{}_{\un{g}'',\un{h}'} \,
     \Gamma^{\un{g}''}{}_{g',ad(\un{h}') \un{h}}       \nonumber \\
   & & + ( C^{\un{g}''}{}_{ad(\un{h}') \un{h}, \un{h}'}
       - C^{\un{g}''}{}_{\un{h}',\un{h}} ) \,
        \Gamma^g{}_{g',\un{g}''} \rbrack \,
        \theta^{\un{h}} \otimes_{\cal A} \theta^{\un{h}'}
\ee
using $\Omega^2 \cong \mbox{im} \, {\bf A}$ (cf section 3) and
(\ref{sigma-theta}).
The condition of vanishing torsion for a linear connection is
\be
     \Gamma^h{}_{g,g'} - \Gamma^h{}_{ad(g) g', g}
   = C^h{}_{g,g'} - C^h{}_{ad(g) g', g}
   = - \delta^h_{g'} + \delta^h_{ad(g) g'}
\ee
which, for a commutative group, reduces to $\Gamma^h{}_{g,g'} =
\Gamma^h{}_{g',g}$.
                                           \hfill {\Large $\Box$}
\vskip.3cm
\noindent
{\em Example 3.}
For a left-covariant first order differential calculus,
\be                              \label{C-connection}
          \Gamma^h{}_{g,g'} = C^h{}_{g,g'}
\ee
defines a left-invariant linear connection which we call the
$C$-{\em connection}. Independent of the continuation of
the first order calculus to higher orders, for this connection the
torsion vanishes.
For a bicovariant differential calculus, the $C$-connection is
bi-invariant as a consequence of (\ref{C-adinv}).
\hspace*{1cm}                                  \hfill {\Large $\Box$}
\vskip.3cm

In Appendix A we introduced the notion of an `extensible connection'.
An extensible linear connection induces a connection on the tensor
product $\Omega^1 \otimes_{\cal A} \Omega^1$ which then enables us
to construct `covariant derivatives' of tensor fields.
In the following we elaborate this notion for the
case of linear connections $\nabla \, : \, \Omega^1 \rightarrow
\Omega^1 \otimes_{\cal A} \Omega^1$ on a finite group $G$ with
a bicovariant (first order) differential calculus (with space of
1-forms $\Omega^1$). Using Proposition \ref{Ext_con}.2, the
following characterization of such connections is obtained.
\vskip.3cm
\noindent
{\em Proposition \ref{Lincon}.1.} A linear connection on a finite group
with a bicovariant (first order) differential calculus is extensible if
and only if there exist bimodule homomorphisms
\be
     V \, : \, \Omega^1 \otimes_{\cal A} \Omega^1 \rightarrow
            \Omega^1 \otimes_{\cal A} \Omega^1 \; , \qquad
     W \, : \, \Omega^1 \rightarrow \Omega^1 \otimes_{\cal A}
            \Omega^1
\ee
such that
\be
  \nabla \varphi = \nabla^\sigma \varphi
                  + V(\varphi \otimes_{\cal A} \rho) + W(\varphi)
\ee
where $\nabla^\sigma$ denotes the connection (\ref{canon-conn}) and
$\rho = \theta^{\un{g}}$.                 \hfill  {\Large $\Box$}
\vskip.3cm
\noindent
It remains to determine the most general form of the bimodule
homomorphisms $V$ and $W$.
\vskip.3cm
\noindent
{\em Proposition \ref{Lincon}.2.} Let $(\Omega^1, \d)$ be a
left-covariant first order differential calculus on a finite group. \\
(a) A map $V \, : \, \Omega^1 \otimes_{\cal A} \Omega^1 \rightarrow
\Omega^1 \otimes_{\cal A} \Omega^1$ is a bimodule homomorphism
if and only if
\be
   V(\theta^g \otimes_{\cal A} \theta^{g'}) = \sum_{h,h'\in \hat{G}
   \atop hh'=g'g} V^{g,g'}_{h,h'} \, \theta^{h'} \otimes_{\cal A}
   \theta^h  \qquad \forall g,g' \in \hat{G}
\ee
with $V^{g,g'}_{h,h'} \in {\cal A}$.    \\
(b) A map $W \, : \, \Omega^1 \rightarrow \Omega^1 \otimes_{\cal A}
\Omega^1$ is a bimodule homomorphism if and only if
\be
   W(\theta^g) = \sum_{h,h'\in \hat{G}
   \atop hh'=g} W^g_{h,h'} \, \theta^{h'} \otimes_{\cal A}
   \theta^h  \qquad \forall g \in \hat{G}
\ee
with $W^g_{h,h'} \in {\cal A}$.
\vskip.1cm
\noindent
{\em Proof:} The proofs of (a) and (b) are essentially the same. We
therefore only present the proof of (b).
Since $\lbrace \theta^g \otimes_{\cal A} \theta^{g'}
\mid g,g' \in \hat{G} \rbrace$ is a left $\cal A$-module basis of
$\Omega^1 \otimes_{\cal A} \Omega^1$, $W(\theta^g)$ must have the form
\begin{eqnarray*}
  W(\theta^g) = \sum_{h,h'\in \hat{G}} W^g_{h,h'} \, \theta^{h'}
                \otimes_{\cal A} \theta^h
\end{eqnarray*}
with coefficients in $\cal A$.
This extends to a left $\cal A$-module homomorphism. The condition for
$W$ to be also right $\cal A$-linear is
\begin{eqnarray*}
 0 = W(\theta^g) \, f - (R_{g^{-1}} f) \, W(\theta^g)
   = \sum_{h,h'\in \hat{G}} (R_{h'^{-1}} R_{h^{-1}} f - R_{g^{-1}} f)
     \; W^g_{h,h'} \, \theta^{h'} \otimes_{\cal A} \theta^h
\end{eqnarray*}
for all $g \in \hat{G}$ and all $f \in {\cal A}$. It is then sufficient
to have this property for all generators $e_{g'}$ (${g'} \in G$) of
$\cal A$, i.e.,
\begin{eqnarray*}
  \sum_{h,h'\in \hat{G}} (e_{g'hh'} - e_{g'g}) \, W^g_{h,h'} \,
  \theta^{h'} \otimes_{\cal A} \theta^h  = 0  \; .
\end{eqnarray*}
This is equivalent to $W^g_{h,h'} =0$ whenever $hh'\neq g$.
                                        \hfill  {\Large $\Box$}
\vskip.3cm
\noindent
According to the two propositions, an extensible linear connection
$\nabla$ is given by
\be
 \nabla \theta^g = \sum_{g',h,h'\in \hat{G} \atop hh'=g'g}
                    V^{g,g'}_{h,h'} \, \theta^{h'} \otimes_{\cal A}
                    \theta^h + \sum_{h,h'\in \hat{G} \atop hh'=g}
                    W^g_{h,h'} \, \theta^{h'} \otimes_{\cal A} \theta^h
\ee
taking into account that $\nabla^\sigma \theta^g = 0$. This means
\be
 \Gamma^g{}_{h,h'} = \left\lbrace \begin{array}{r@{\quad}c@{\quad}l}
- W^g_{h,h'}             &             & hh' = g  \\
- V^{g,hh'g^{-1}}_{h,h'} &  \mbox{for} & g' = hh'g^{-1} \in \hat{G} \\
0                        &             & hh'g^{-1} \not\in \hat{G} \cup
                                         \{ e \}
                                   \end{array} \right.
\ee
from which we observe that there are restrictions for $\nabla$ to
be extensible iff there are products $hh'g^{-1} \not\in \hat{G} \cup
\{ e \}$ for $h,h',g \in \hat{G}$.
\vskip.2cm
\noindent
{\em Example 4.} For $G = \Ir_4 = \lbrace e, a, a^2, a^3 \rbrace$ and
$\hat{G} = \lbrace a, a^2 \rbrace$ we find the restrictions
\be
     \Gamma^{a^2}_{a,a} = \Gamma^a_{a^2,a^2} = 0
\ee
for a linear connection to be extensible. This excludes the
$C$-connection of example 3.                  \hfill {\Large $\Box$}

\section{Vector fields, dual connections, and metrics on finite groups}
\label{Vecfield}
\setcounter{equation}{0}
Let $\Omega$ be a left-covariant differential calculus on a finite
group $G$. By $\cal X$ we denote the dual $\cal A$-bimodule of
$\Omega^1$ with duality contraction $\langle \varphi , X \rangle$ for
$\varphi \in \Omega^1$ and $X \in {\cal X}$ (see Appendix B).
The elements of $\cal X$ act as operators on $\cal A$ via
\be
    X f := \langle \d f , X \rangle  \; .
\ee
The (nonvanishing) Maurer-Cartan forms $\theta^g$ constitute a basis
of $\Omega^1$ as a left or right $\cal A$-module. Let
$\lbrace \ell'_g \mid g \in \hat{G} \rbrace$ be the dual basis. Then
\be
     \ell'_g f
   = \langle \d f , \ell'_g \rangle
   = \langle (\ell_{\un{h}} f) \, \theta^{\un{h}} , \ell'_g \rangle
   = \ell_g f
\ee
shows that $\ell'_g = \ell_g$.
In the same way one verifies that $\lbrace r_g \mid g \in \hat{G}
\rbrace$ is the dual basis of $\lbrace \omega^g \rbrace$.
Elements of ${\cal X}$ can now be written as
\be
   X = \ell_{\un{g}} \cdot X^{\un{g}}
\ee
where  $X f = (\ell_{\un{g}} f) \, X^{\un{g}}$. As a consequence of
(\ref{f-theta}),
\be
 \langle \theta^h , ({\cal R}_g f) \, \ell_g - \ell_g \cdot f \rangle
 = \langle \theta^h ({\cal R}_g f) , \ell_g \rangle
   - \delta^h_g \, f
 = (R_{h^{-1} g} f) \langle \theta^h , \ell_g \rangle
   - \delta^h_g \, f
 = 0
\ee
so that
\be
    ({\cal R}_g f) \, \ell_g = \ell_g \cdot f  \; .
\ee
\vskip.3cm

By duality (see Appendix B) a linear left $\cal A$-module connection
$\nabla$ induces a right $\cal A$-module connection $\nabla^\ast$ on
$\cal X$ such that
\be
    \nabla^\ast \ell_g = \ell_{\un{h}} \otimes_{\cal A}
    \omega^{\un{h}}{}_g
\ee
where $\omega^g{}_{g'}$ are the connection 1-forms with respect to
the basis $\theta^g$ (cf (\ref{conn-1forms})). $\nabla^\ast$ extends
to a map ${\cal X} \otimes_{\cal A} \Omega \rightarrow {\cal X}
\otimes_{\cal A} \Omega$ such that
\be
   \nabla^\ast (\chi \varphi) = (\nabla^\ast \chi) \, \varphi
                               + \chi \, \d \varphi
\ee
for $\chi \in {\cal X} \otimes_{\cal A} \Omega$ and $\varphi \in
\Omega$. Regarding\footnote{We may also regard $\Xi$ as an element
of the product module ${\cal X} \otimes_{\cal A} \Omega^1$. A
connection $\nabla$ on $\Omega^1$ together with its dual $\nabla^\ast$
can then be used to define $\tilde{\nabla}(X \otimes_{\cal A} \varphi)
:= (\nabla^\ast X) \otimes_{\cal A} \varphi + X \otimes_{\cal A} \nabla
\varphi$ which implies $\tilde{\nabla} \Xi = 0$. This makes sense
though $\tilde{\nabla}$ is {\em not} a left or right $\cal A$-module
connection on ${\cal X} \otimes_{\cal A} \Omega^1$. }
the {\em canonical form}
\be
   \Xi := \ell_{\un{g}} \otimes_{\cal A} \theta^{\un{g}}
\ee
as an element of ${\cal X} \otimes_{\cal A} \Omega$, we find
\be
   \nabla^\ast \Xi = \nabla^\ast \ell_{\un{g}} \otimes_{\cal A}
   \theta^{\un{g}} + \ell_{\un{g}} \otimes_{\cal A} \d \theta^{\un{g}}
   = \ell_{\un{g}} \otimes_{\cal A} D \theta^{\un{g}}
\ee
with
\be
   D \theta^g := \d \theta^g + \omega^g{}_{\un{h}} \, \theta^{\un{h}}
               =: \Theta^g    \; .
\ee
This is another (equivalent) expression for the torsion of $\nabla$.
Furthermore, since $\nabla^\ast \Xi$ is again an element of
${\cal X} \otimes_{\cal A} \Omega$, we can apply $\nabla^\ast$
another time. This yields
\be
 (\nabla^\ast)^2 \, \Xi = \ell_{\un{g}} \otimes_{\cal A}
                          D \Theta^{\un{g}}
\ee
where
\be
 D \Theta^g = \d \Theta^g + \omega^g{}_{\un{g}'} \, \Theta^{\un{g}'}
            = D^2 \theta^g = \Omega^g{}_{\un{g}'} \, \theta^{\un{g}'}
\ee
which resembles the first Bianchi identity of classical differential
geometry.
\vskip.2cm

As a {\em metric} we may regard an element ${\bf g} \in{\cal X}
\otimes_{\cal A}{\cal X}$ with certain properties.\footnote{A reality
or hermiticity condition requires an involution, an extra structure
which we leave aside in the present work. Less straightforward is
the implementation of a notion of invertibility. Note that we could
also think of a metric as an element of $\Omega^1 \otimes_{\cal A}
\Omega^1$. Then Appendix \ref{Invtensor} provides us with examples.
See also the discussion in section \ref{Concl}.}
In terms of the basis $\ell_g \otimes_{\cal A} \ell_{g'}$ we have
${\bf g} = \ell_{\un{g}} \otimes_{\cal A} \ell_{\un{g}'} \cdot
{\bf g}^{\un{g},\un{g}'}$.
A metric is called {\em compatible} with a connection on
${\cal X} \otimes_{\cal A} {\cal X}$ if $\bf g$ is covariantly constant.
\vskip.3cm
\noindent
{\em Example.}
For a bicovariant differential calculus, the canonical bimodule
isomorphism $\sigma$ has a `dual' $\sigma' \, : \, \Omega^1
\otimes_{\cal A} {\cal X} \rightarrow {\cal X} \otimes_{\cal A}
\Omega^1$ (cf (\ref{dual-Psi})). From
\be
  \langle \theta^{g'} , \sigma' (\theta^h \otimes_{\cal A} \ell_g)
  \rangle
  &=& \langle \sigma(\theta^{g'} \otimes_{\cal A} \theta^h) , \ell_g
    \rangle
  = \langle \theta^{g'^{-1}hg'} \otimes_{\cal A} \theta^{g'} , \ell_g
    \rangle                                        \nonumber \\
  &=& \theta^{g'^{-1}hg'} \, \delta^{g'}_g
  = \langle \theta^{g'} , \ell_g \otimes_{\cal A} \theta^{g^{-1}hg}
    \rangle
\ee
we deduce
\be
    \sigma'(\theta^h \otimes_{\cal A} \ell_g)
  = \ell_g \otimes_{\cal A} \theta^{g^{-1} h g}  \; .
\ee
For the connection $\nabla^\sigma$ defined in (\ref{canon-conn})
the dual connection (see Appendix B) is given by
\be
  \nabla^{\sigma'} X = X \otimes_{\cal A} \rho
                       - \sigma'(\rho \otimes_{\cal A} X)
\ee
and has the property $\nabla^{\sigma'} \ell_g = 0$. It is extensible
and we obtain
\be
   \nabla^{\sigma'}_\otimes {\bf g} = \ell_{\un{g}} \otimes_{\cal A}
   \ell_{\un{g}'} \cdot \d {\bf g}^{\un{g},\un{g}'}
\ee
so that $\bf g$ is compatible with $\nabla^{\sigma'}_\otimes$ if
and only if $\d {\bf g}^{\un{g},\un{g}'} = 0$.\footnote{The equation
$\d f = 0$ for $f \in {\cal A}$ does not necessarily imply $f \in \Cx$.
It depends on the differential calculus what the `constant functions'
are.}

Alternatively, we can extend $\nabla^\sigma$ to a connection
$\nabla^\sigma_\otimes$ on $\Omega^1 \otimes_{\cal A} \Omega^1$ and
there is a dual $\nabla^\sigma_\otimes{}^\ast$ of the latter
(see Appendix B). Then
\be
 \langle \theta^g \otimes_{\cal A} \theta^{g'} ,
 \nabla^\sigma_\otimes{}^\ast {\bf g} \rangle = \d {\bf g}^{g,g'}
  - \langle \nabla^\sigma_\otimes (\theta^g \otimes_{\cal A}
  \theta^{g'}) , {\bf g} \rangle
\ee
for a metric $\bf g$. Since $\theta^g$ is covariantly constant with
respect to $\nabla^\sigma$, we again obtain $\d {\bf g}^{g,g'} = 0$
as the metric compatibility condition.
                                                \hfill {\Large $\Box$}
\vskip.3cm

Let $\Omega^1$ be the space of 1-forms of a bicovariant first order
differential calculus on a finite group.
Symmetry conditions can then be imposed on a metric as follows. For
example, we call $\bf g$ s-{\em symmetric} if $\langle \sigma
(\varphi \otimes_{\cal A} \varphi') , {\bf g} \rangle = \langle \varphi
\otimes_{\cal A} \varphi' , {\bf g} \rangle$ for all $\varphi ,
\varphi' \in \Omega^1$. This can also be expressed as $\sigma_{\cal X}
{\bf g} = {\bf g}$ in terms of the transpose $\sigma_{\cal X}$ of
$\sigma$, which is determined by
\be
   \langle \theta^g \otimes_{\cal A} \theta^{g'} , \sigma_{\cal X}
      ( \ell_h \otimes_{\cal A} \ell_{h'} ) \rangle
  &=& \langle \sigma(\theta^g \otimes_{\cal A} \theta^{g'}) ,
      \ell_h \otimes_{\cal A} \ell_{h'}  \rangle   \nonumber \\
  &=& \langle \theta^{ad(g^{-1}) g'} \otimes_{\cal A} \theta^g ,
      \ell_h \otimes_{\cal A} \ell_{h'}  \rangle   \nonumber \\
  &=& \delta^g_h \, \delta^{ad(g^{-1}) g'}_{h'}
   =  \delta^g_h \, \delta_{ad(h) h'}^{g'}         \nonumber \\
  &=& \langle \theta^g \otimes_{\cal A} \theta^{g'} ,
      \ell_{ad(h)h'} \otimes_{\cal A} \ell_h  \rangle
\ee
i.e.,
\be
    \sigma_{\cal X} ( \ell_g \otimes_{\cal A} \ell_{g'} )
  = \ell_{ad(g)g'} \otimes_{\cal A} \ell_g   \; .
\ee
Using $r_g = \ell_{ad(\un{h}^{-1}) g} \cdot e_{\un{h}}$ and
$\ell_g \cdot e_h = e_{hg^{-1}} \, \ell_g$ one finds
$\sigma_{\cal X} ( \ell_g \otimes_{\cal A} r_{g'} )
 = r_{g'} \otimes_{\cal A} \ell_g$ which is the analogue of
(\ref{sig_th_om}). Indeed, $\cal X$ inherits from $\Omega^1$ the
structure of a bicovariant bimodule and $\sigma_{\cal X}$ is the
corresponding braid operator.
\vskip.3cm

According to the general construction in Appendix C, a left (right)
coaction on $\Omega^1$ induces a left (right) coaction on $\cal X$.
For a left-covariant first order differential calculus the left
coaction on $\cal X$ is determined by
\be
        \Delta_{\cal X} (\ell_g) = \idty \otimes \ell_g
\ee
i.e., the elements dual to the left-invariant basis 1-forms $\theta^g$
are also left-invariant. An element $X = \ell_{\un{g}} \cdot
X^{\un{g}}$ of $\cal X$ is left-invariant iff $X^g \in \Cx$ for all
$g \in \hat{G}$. The coaction extends to ${\cal X} \otimes_{\cal A}
{\cal X}$ (see Appendix C). A tensor ${\bf g} = \ell_{\un{g}}
\otimes_{\cal A} \ell_{\un{g}'} \cdot {\bf g}^{\un{g},\un{g}'}$ is
then left-invariant iff ${\bf g}^{\un{g},\un{g}'} \in \Cx$.
The above example now shows that every left-invariant metric is
covariantly constant with respect to (the extension of) the connection
$\nabla^{\sigma'}$.

\section{Noncommutative geometry of the symmetric group ${\cal S}_3$}
\label{NgS3}
\setcounter{equation}{0}
We denote the elements of the symmetric group ${\cal S}_3$ as follows,
\be
        a = (12)  \, , \quad
        b = (23)  \, , \quad
        c = (13)  \, ,
\ee
$ab, ba$, and $e$ (the unit element). In order to determine the
${\cal S}_3$ left-invariant differential calculi on a set of six
elements, we have to calculate the orbits of the left-action on
$({\cal S}_3 \times {\cal S}_3)'$. These are
\begin{eqnarray*}
   {\cal O}_1 &=& \lbrace (e,a), (a,e), (b,ba), (c,ab), (ab,c),
                          (ba,b) \rbrace       \\
   {\cal O}_2 &=& \lbrace (e,b), (a,ab), (b,e), (c,ba), (ab,a),
                          (ba,c) \rbrace       \\
   {\cal O}_3 &=& \lbrace (e,c), (a,ba), (b,ab), (c,e), (ab,b),
                          (ba,a) \rbrace       \\
   {\cal O}_4 &=& \lbrace (e,ab), (a,b), (b,c), (c,a), (ab,ba),
                          (ba,e) \rbrace       \\
   {\cal O}_5 &=& \lbrace (e,ba), (a,c), (b,a), (c,b), (ab,e),
                          (ba,ab) \rbrace
\end{eqnarray*}
They are in correspondence with elements of $G \setminus \lbrace e
\rbrace$. Left-invariant differential calculi are obtained by deleting
subgraphs corresponding to orbits from the graph
corresponding to the universal differential calculus (which is
left-invariant, of course).
\vskip.2cm

With respect to left- and right-action $({\cal S}_3 \times {\cal S}_3)'$
decomposes into two orbits,\footnote{These orbits are in correspondence
with the nontrivial conjugacy classes in ${\cal S}_3$, i.e.,
$\lbrace a,b,c \rbrace$ and $\lbrace ab, ba \rbrace$ (cf section
\ref{Diffcalc}).}
\be
   {\cal O}_I = {\cal O}_1 \cup {\cal O}_2 \cup {\cal O}_3 \, ,
   \quad
   {\cal O}_{II} = {\cal O}_4 \cup {\cal O}_5   \; .
\ee
Hence, there are two bicovariant differential calculi, $\Omega^1_{I}$
and $\Omega^1_{II}$, on ${\cal S}_3$ besides the universal and the
trivial one. Their graphs are obtained by deleting all arrows
corresponding to elements of ${\cal O}_I$ or ${\cal O}_{II}$,
respectively, from the graph associated with the universal differential
calculus (see Fig. 2). All these graphs are symmetric in the sense that
for every arrow also the reverse arrow is present.

\unitlength1.cm
\begin{picture}(12.,3.2)(-1.4,-0.3)
\thicklines
\put(1.,1.) {\circle*{0.2}}
\put(2.,0.) {\circle*{0.2}}
\put(2.,2.) {\circle*{0.2}}
\put(3.5,0.) {\circle*{0.2}}
\put(3.5,2.) {\circle*{0.2}}
\put(4.5,1.) {\circle*{0.2}}
\put(0.6,0.9) {$e$}
\put(1.9,-0.4) {$a$}
\put(1.9,2.3) {$ab$}
\put(3.4,-0.4) {$b$}
\put(3.4,2.3) {$ba$}
\put(4.8,0.9) {$c$}
\put(1.,1.) {\line(1,1){1.}}
\put(1.,1.) {\line(5,2){2.5}}
\put(2.,2.) {\line(1,0){1.5}}
\put(2.,0.) {\line(1,0){1.5}}
\put(2.,0.) {\line(5,2){2.5}}
\put(3.5,0.) {\line(1,1){1.}}
%%%
\put(7.,1.) {\circle*{0.2}}
\put(8.,0.) {\circle*{0.2}}
\put(8.,2.) {\circle*{0.2}}
\put(9.5,0.) {\circle*{0.2}}
\put(9.5,2.) {\circle*{0.2}}
\put(10.5,1.) {\circle*{0.2}}
\put(6.6,0.9) {$e$}
\put(7.9,-0.4) {$a$}
\put(7.9,2.3) {$ab$}
\put(9.4,-0.4) {$b$}
\put(9.4,2.3) {$ba$}
\put(10.8,0.9) {$c$}
\put(7.,1.) {\line(1,-1){1.}}
\put(7.,1.) {\line(1,0){3.5}}
\put(7.,1.) {\line(5,-2){2.5}}
\put(8.,0.) {\line(0,1){2.}}
\put(8.,0.) {\line(3,4){1.5}}
\put(9.5,0.) {\line(-3,4){1.5}}
\put(9.5,0.) {\line(0,1){2.}}
\put(10.5,1.) {\line(-5,2){2.5}}
\put(10.5,1.) {\line(-1,1){1.}}
\end{picture}

\begin{center}
\begin{minipage}{10cm}
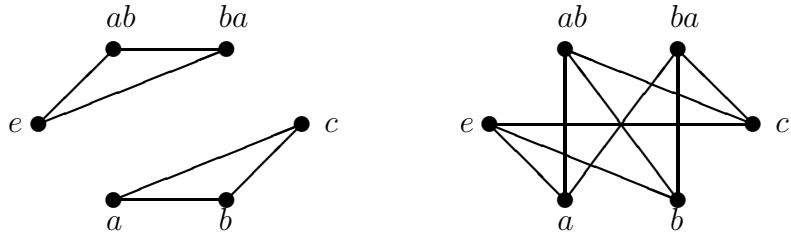
\centerline{\bf Fig. 2}
\vskip.1cm
\noindent
The graphs which determine the bicovariant first order differential
calculi $\Omega^1_I$ and $\Omega^1_{II}$ on ${\cal S}_3$.
\end{minipage}
\end{center}
\vskip.3cm
\noindent
The bimodule isomorphism $\sigma$ is a nontrivial map in the case
under consideration. In terms of the decomposition into conjugacy
classes
\be
   {\cal S}_3 = \lbrace e \rbrace \cup \underbrace{\lbrace a,b,c
   \rbrace}_{\mbox{$=:{\cal S}_3'$}} \cup \underbrace{\lbrace ab,ba
   \rbrace}_{\mbox{$=:{\cal S}_3''$}}
\ee
it is given by
\be \label{sigs3}
  \sigma(\theta^x, \theta^x) &=& \theta^x \otimes_{\cal A} \theta^x
       \qquad  \forall x \in {\cal S}_3     \nonumber \\
  \sigma(\theta^x, \theta^y) &=& \theta^z \otimes_{\cal A} \theta^x
       \qquad  \forall x,y \in {\cal S}_3', x \neq y, z \in {\cal S}_3'
  \setminus \lbrace x,y \rbrace           \nonumber \\
  \sigma(\theta^{ab}, \theta^{ba}) &=& \theta^{ba} \otimes_{\cal A}
       \theta^{ab}        \nonumber \\
  \sigma(\theta^{ba}, \theta^{ab}) &=& \theta^{ab} \otimes_{\cal A}
       \theta^{ba}        \nonumber \\
  \sigma(\theta^x, \theta^{ab}) &=& \theta^{ba} \otimes_{\cal A}
       \theta^x   \qquad  \forall x \in {\cal S}_3'  \nonumber \\
  \sigma(\theta^x, \theta^{ba}) &=& \theta^{ab} \otimes_{\cal A}
       \theta^x   \qquad  \forall x \in {\cal S}_3'  \nonumber \\
  \sigma(\theta^{ab}, \theta^x) &=& \theta^y \otimes_{\cal A}
       \theta^{ab}  \qquad \mbox{for } (x,y) \in \lbrace (a,c), (b,a),
       (c,b) \rbrace                                     \nonumber \\
  \sigma(\theta^{ba}, \theta^x) &=& \theta^y \otimes_{\cal A}
       \theta^{ba}  \qquad \mbox{for } (x,y) \in \lbrace (a,b), (b,c)
       (c,a) \rbrace
\ee
for the universal first order differential calculus. Since the center
of ${\cal S}_3$ is trivial, $ad({\cal S}_3) \cong {\cal S}_3$ and
$|ad({\cal S}_3)|=6$, so that $\sigma^{12} = id$ according to
(\ref{sigid}). For the other bicovariant calculi, the corresponding
$\sigma$ is induced in an obvious way. In the case of $\Omega^1_{I}$
the bimodule isomorphism is the one of the commutative (sub)group
$\Ir_3$. Hence, $\sigma_{I}^2=id$. For $\Omega^1_{II}$ one can deduce
from (\ref{sigs3}) that $\sigma_{II}^3 = id$. For the restriction of
$\sigma$ to the subbimodule of $\tilde{\Omega}^1 \otimes_{\cal A}
\tilde{\Omega}^1$ which is generated by $\lbrace \theta^x
\otimes_{\cal A} \theta^{x'}, \theta^{x'} \otimes_{\cal A} \theta^{x}
\mid x \in {\cal S}_3', x' \in {\cal S}_3'' \rbrace$ we have
$\sigma^4=id$. Hence $2 \, |ad({\cal S}_3)| = 12$ is actually the order
of $\sigma$ for the universal first order differential calculus, in
accordance with Proposition 3.2.
\vskip.2cm

In order to determine the most general bi-invariant linear connection
on ${\cal S}_3$ with the universal first order differential calculus,
one has to determine the $(ad)^3$-orbits in $(G \setminus \lbrace e
\rbrace)^3$. There are 24 of them so that bi-invariance restricts the
125 connection coefficients in (\ref{lin_conn}) to 24 independent
constants.

\subsection{Geometry of the 3-dimensional bicovariant calculus}
The bimodule $\Omega^1_{II}$ is generated (as a left $\cal A$-module)
by $\theta^a, \theta^b, \theta^c$. The operator $\sigma$ acts on
$\Omega^1_{II} \otimes_{\cal A} \Omega^1_{II}$ as follows,
\be
      \sigma(\theta^x \otimes_{\cal A} \theta^x)
  &=& \theta^x \otimes_{\cal A} \theta^x
      \qquad  \forall x \in \lbrace a,b,c \rbrace  \nonumber \\
      \sigma(\theta^x \otimes_{\cal A} \theta^y)
  &=& \theta^z \otimes_{\cal A} \theta^x
\ee
where, in the last equation, $z$ is the complement of $x,y$ in
${\cal S}'_3$. As already mentioned, $\sigma^2 \neq id$, but
$\sigma^3 = id$. This rules out $-1$ as an eigenvalue of
$\sigma$ so that all elements of $\Omega^1_{II} \otimes_{\cal A}
\Omega^1_{II}$ are w-symmetric according to Proposition 3.3.
The bimodule $\Omega^1_{II} \otimes_{\cal A} \Omega^1_{II}$ splits
into a direct sum of subbimodules (cf Proposition 3.3),
$ \Omega^1_{II} \otimes_{\cal A} \Omega^1_{II} = \mbox{ker} {\bf A}
\oplus \mbox{im} {\bf A}$ where $\mbox{ker} {\bf A}$ is generated by
$\theta^x \otimes_{\cal A}
\theta^x$ for $x \in \lbrace a,b,c \rbrace$, $\theta^a \otimes_{\cal A}
\theta^b + \theta^b \otimes_{\cal A} \theta^c + \theta^c
\otimes_{\cal A} \theta^a$ and $\theta^b \otimes_{\cal A}
\theta^a + \theta^a \otimes_{\cal A} \theta^c + \theta^c
\otimes_{\cal A} \theta^b$. These are eigenvectors of $\sigma$ with
eigenvalue $1$, i.e., s-symmetric tensors.
The image of $\Omega^1_{II} \otimes_{\cal A} \Omega^1_{II}$ under
$\bf A$ is generated by $\theta^a \otimes_{\cal A} \theta^b
- \theta^c \otimes_{\cal A} \theta^a$, $\theta^a \otimes_{\cal A}
\theta^b - \theta^b \otimes_{\cal A} \theta^c$, $\theta^a
\otimes_{\cal A} \theta^c - \theta^b \otimes_{\cal A} \theta^a$ and
$\theta^a \otimes_{\cal A} \theta^c - \theta^c \otimes_{\cal A}
\theta^b$. These are w-antisymmetric tensors. The space of 2-forms
({\`a} la Woronowicz) is therefore four-dimensional. A basis is given
by $\theta^a \theta^b, \, \theta^b \theta^c, \, \theta^a \theta^c, \,
\theta^b \theta^a$ and we have the relations
\be
   \theta^c \theta^a = - \theta^a \theta^b - \theta^b \theta^c
               \; , \qquad
   \theta^c \theta^b = - \theta^b \theta^a - \theta^a \theta^c
               \; , \qquad
   \theta^x \theta^x = 0 \qquad x \in \lbrace a,b,c \rbrace  \; .
\ee
\vskip.2cm

For a linear connection, there are a priori 27 connection coefficients.
Bi-invariance restricts them as follows,
\be
 \Gamma^a_{a,a} &=& \Gamma^b_{b,b} = \Gamma^c_{c,c}   \nonumber \\
 \Gamma^a_{a,b} &=& \Gamma^a_{a,c} = \Gamma^b_{b,a} = \Gamma^b_{b,c}
                 = \Gamma^c_{c,a} = \Gamma^c_{c,b}  \nonumber \\
 \Gamma^a_{b,a} &=& \Gamma^a_{c,a} = \Gamma^b_{a,b} = \Gamma^b_{c,b}
                 = \Gamma^c_{a,c} = \Gamma^c_{b,c}  \nonumber \\
 \Gamma^a_{b,b} &=& \Gamma^a_{c,c} = \Gamma^b_{a,a} = \Gamma^b_{c,c}
                 = \Gamma^c_{a,a} = \Gamma^c_{b,b}  \nonumber \\
 \Gamma^a_{b,c} &=& \Gamma^a_{c,b} = \Gamma^b_{a,c} = \Gamma^b_{c,a}
                 = \Gamma^c_{a,b} = \Gamma^c_{b,a}   \; .
\ee
The condition of vanishing torsion for a bi-invariant connection becomes
\be
   \Gamma^a_{b,a} = \Gamma^a_{a,b} = \Gamma^a_{b,c} - 1
\ee
leaving us with only three independent constants. It turns out that
$\Gamma^x{}_{y,z} = \Gamma^x{}_{z,y}$ for bi-invariant connections
without torsion. Among these is the $C$-connection
(\ref{C-connection}) for which $\Gamma^a_{a,a} = -2, \, \Gamma^a_{a,b}
= -1, \, \Gamma^a_{b,b} = \Gamma^a_{b,c} = 0$. There are no
bi-invariant connections for which torsion and curvature vanish.
\vskip.2cm

In the case under consideration, we have $\hat{G} = \lbrace a,b,c
\rbrace = \hat{G}^{-1}$, $\hat{G} \cdot \hat{G} = \lbrace e,ab,ba
\rbrace$ and $\hat{G} \cdot \hat{G} \cdot \hat{G} = \hat{G}$. As a
consequence, every linear connection is extensible (cf section
\ref{Lincon}).
But there are no (non-trivial) bimodule homomorphisms $W \, : \,
\Omega^1_{II} \rightarrow \Omega^1_{II} \otimes_{\cal A} \Omega^1_{II}$
since $gg' \not\in \hat{G}$ for all $g,g' \in \hat{G}$. From Proposition
\ref{Lincon}.1 we infer that all linear connections have the form
\be
   \nabla \varphi = \rho \otimes_{\cal A} \varphi
                   + V(\varphi \otimes_{\cal A} \rho)
\ee
with a bimodule homomorphism $V$.
This includes the family (\ref{conn-fam}), of course, which in the
case under consideration depends on three independent constants.
The only linear (left $\cal A$-module) connection with right
$\cal A$-linearity is given by $\nabla \varphi = \rho \otimes_{\cal A}
\varphi$ in accordance with Proposition \ref{Ext_con}.3.

\section{Final remarks}
\label{Concl}
\setcounter{equation}{0}
In this work we have continued our previous research on
noncommutative geometry of discrete sets \cite{DMH94-graphs,DMHV}
and, in particular, finite groups \cite{Sita94-groups,DMH94-groups}.
\vskip.2cm

Much of the material presented concerns the notion of linear
connections.
In \cite{Mour95} a special class of linear left $\cal A$-module
connections has been considered satisfying the additional condition
\be                                   \label{Mour-cond}
   \nabla(\varphi \, f) = \tau (\varphi \otimes_{\cal A} \d f)
                          + (\nabla \varphi) \, f
\ee
where $\tau$ is an $\cal A$-bimodule homomorphism. In classical
differential geometry, all linear connections satisfy this condition
with the choice of the permutation map for $\tau$. This observation
was taken in \cite{Mour95} to consider the above condition in
noncommutative geometry. It should be noticed, however, that
connections in commutative geometry {\em automatically} satisfy
this condition whereas in noncommutative geometry it severely
restricts the possible (linear) connections, in general. It is
therefore quite unclear from this point of view what the relevance of
the class of (linear) connections determined by (\ref{Mour-cond}) is.
However, it has also been pointed out in \cite{Mour95} that
linear connections with the above property can be extended
to tensor products (over $\cal A$) of 1-forms. Indeed, given connections
on two $\cal A$-bimodules, it seems to be impossible, in general,
to build from these a connection on the tensor product of the
two modules. In particular, we would like to achieve this in order
to be able to talk about a covariantly constant metric. In Appendix A
we have addressed the question of extensibility in more generality.
In our attempt to solve this problem, we were led to the condition
(\ref{Mour-cond}), which we therefore called `extensibility
condition'. This provides a much stronger motivation for the
consideration of the special class of linear connections satisfying
(\ref{Mour-cond}). We have to stress, however, that there may still
be a way beyond our ansatz to extend connections. In Appendix D we
briefly discussed a natural modification of the usual definition
of a linear connection which guarantees extensibility. It turned
out to be too restrictive, however.
\vskip.2cm

For a bicovariant differential calculus on a Hopf algebra there is
a canonical choice for $\tau$, the canonical bimodule isomorphism
$\sigma$ \cite{Woro89}. Using the fact that powers of $\sigma$ are
again bimodule isomorphisms, one actually has a whole class of
extensible linear connections on the Hopf algebra. Similar
observations have been made in \cite{GMMM95} where, however,
the restriction to `generalized permutations' $\tau$ which satisfy
$\pi \circ (\tau + id) = 0$ rules out $\tau = \sigma$ (together with
$\pi = (id - \sigma)/2$ which is used in \cite{Woro89} to extend
bicovariant first order differential calculi on Hopf algebras
to higher orders) if $\sigma^2 \neq id$.
The reasoning behind this restriction (see also \cite{DVMMM}) is not
quite transparent for us and in the formalism presented in this paper
(which extends beyond finite groups) there is no natural place for it.
\vskip.2cm

We should stress that extensibility conditions for connections not
only arise for noncommutative algebras, but already for commutative
algebras with `noncommutative differential calculi' (where functions
do not commute with 1-forms, in general), hence in particular for
differential calculi on finite sets. For finite groups we have
elaborated the extension condition for linear connections and worked
out the corresponding restrictions.
\vskip.2cm

The noncommutativity of a differential calculus with space of 1-forms
$\Omega^1$ results in a non-locality of the tensor product of 1-forms.
This manifests itself in the fact that components of an element
$\alpha \in \Omega^1 \otimes_{\cal A} \Omega^1$ with respect to
some (left or right $\cal A$-module) basis of $\Omega^1$ do not
transform in a covariant manner under a change of basis.\footnote{This
is in contrast to the fact that other basic constructions in
noncommutative geometry indeed lead to quantities with covariant
components, see section \ref{Lincon}.}
Though from a mathematical point of view one can hardly think of an
alternative of the tensor product over the algebra $\cal A$,
it may not be (directly) suitable for a description of physics.
It seems that some modification is needed. An example is
provided by \cite{DMH92-grav} where, for a certain noncommutative
differential calculus on manifolds, a modified wedge product was
constructed with the help of a linear connection, which then allowed to
read off covariant components from (generalized) differential forms.
\vskip.2cm

It is the problem just mentioned which indicates that
in noncommutative geometry the concept of a metric as an element
of $\Omega^1 \otimes_{\cal A} \Omega^1$, respectively a dual module,
may be too naive. A departure from this concept might have crucial
consequences for the relevance of the class of extensible linear
connections, of course. Further exploration of (finite) examples, and
perhaps even those presented in this work, should shed more light on
these problems.
\vskip.2cm

As we have seen in section \ref{Lincon}, there are examples in the
class of extensible connections which have a geometrical meaning as
analogues of the $(\pm)$-parallelisms on Lie groups. More generally,
corresponding connections exist on every Hopf algebra with a bicovariant
differential calculus which is inner with a bi-invariant 1-form
$\rho$, so in particular on the quantum groups $GL_q(n)$ (see also
\cite{GMMM95}). On the other hand, naturally associated with a
left-covariant differential calculus on a finite group is the
$C$-connection (introduced in section 4) which is not extensible,
in general.
\vskip.2cm

We have developed `differential geometry' on finite sets to a level
which now enables us to write down `geometric equations' on discrete
differential manifolds and to look for exact solutions. Comparatively
simple examples are given by the equations of vanishing curvature or
vanishing torsion for a linear connection, in which cases we presented
exact solutions. More interesting would be an analogue of
Einstein's equations, of course. But, as mentioned above, there is
still something to be understood concerning the concept of a metric
before we can seriously proceed towards this goal.
\vskip.3cm
\noindent
{\bf Acknowledgment.} F. M.-H. is grateful to Professor H. L.
deVries for a helpful discussion on symmetric groups.

\begin{appendix}
\renewcommand{\theequation} {\Alph{section}.\arabic{equation}}

\section{Extension of connections to tensor products of bimodules}
\label{Ext_con}
\setcounter{equation}{0}
Let $\cal A$ be an associative unital algebra, $\Gamma$ an
$\cal A$-bimodule, $\Gamma'$ a left $\cal A$-module, and
$\nabla, \nabla'$ left $\cal A$-module connections on $\Gamma$ and
$\Gamma'$, respectively, with respect to a first order differential
calculus on $\cal A$ with space of 1-forms $\Omega^1$. We would like
to build from these a connection on $\Gamma \otimes_{\cal A} \Gamma'$,
i.e., a map
\be
   \nabla \, : \, \Gamma \otimes_{\cal A} \Gamma' \; \rightarrow  \;
   \Omega^1 \otimes_{\cal A} \Gamma \otimes_{\cal A} \Gamma'
\ee
which is $\Cx$-linear and satisfies
\be
   \nabla(f (\gamma \otimes_{\cal A} \gamma')) = \d f \otimes_{\cal A}
   \gamma \otimes_{\cal A} \gamma' + f \, \nabla(\gamma \otimes_{\cal A}
   \gamma')  \qquad \forall f \in {\cal A}, \, \gamma \in \Gamma,
   \, \gamma' \in \Gamma'     \; .         \label{tp-con1}
\ee
In order to be well-defined, the extended connection has to satisfy
\be                         \label{tp-con2}
   \nabla(\gamma f \otimes_{\cal A} \gamma')  =
   \nabla(\gamma \otimes_{\cal A} f \gamma')    \; .
\ee
Let us consider the ansatz
\be                               \label{ext-ansatz}
   \nabla = \Phi \circ (\nabla \otimes id_{\Gamma'}) +
   \Psi \circ (id_{\Gamma} \otimes \nabla')
\ee
with linear maps
\begin{eqnarray*}
   \Phi \, : \, \Omega^1 \otimes_{\cal A} \Gamma
             \otimes_{\cal A} \Gamma' \; & \rightarrow & \;
             \Omega^1 \otimes_{\cal A} \Gamma
             \otimes_{\cal A} \Gamma'                    \\
   \Psi \, : \, \Gamma \otimes_{\cal A} \Omega^1
             \otimes_{\cal A} \Gamma' \; & \rightarrow & \;
             \Omega^1 \otimes_{\cal A} \Gamma
             \otimes_{\cal A} \Gamma'    \; .
\end{eqnarray*}
In the following we evaluate the conditions which the extended
connection has to satisfy. These restrict the possibilities
for the maps $\Phi$ and $\Psi$.
{}From (\ref{tp-con2}) we obtain the following condition,
\be
   \Psi(\gamma \otimes_{\cal A} \d f \otimes_{\cal A} \gamma')
  = \Phi( \lbrack \nabla(\gamma f) - (\nabla \gamma) f \rbrack
    \otimes_{\cal A} \gamma')  \; .      \label{tp-con2+}
\ee
(\ref{tp-con1}) leads to
\be
 0 &=& \Phi( \d f \otimes_{\cal A} \gamma  \otimes_{\cal A} \gamma')
       - \d f \otimes_{\cal A} \gamma \otimes_{\cal A} \gamma'
       \nonumber \\
   & & + \Phi(f \, \nabla \gamma \otimes_{\cal A} \gamma')
       - f \, \Phi(\nabla \gamma \otimes_{\cal A} \gamma')
       \nonumber \\
   & & + \Psi(f \gamma \otimes_{\cal A} \nabla' \gamma')
       - f \, \Psi(\gamma \otimes_{\cal A} \nabla' \gamma')  \; .
\ee
If we demand $\Phi$ and $\Psi$ to be left $\cal A$-linear, then
the last equation implies
\be
  \Phi = id_{\Omega^1} \otimes id_{\Gamma} \otimes id_{\Gamma'}
\ee
and (\ref{tp-con2+}) reduces to
\be
   \Psi = \Psi_\nabla \otimes id_{\Gamma'}
\ee
with a map $\Psi_\nabla \, : \, \Gamma \otimes_{\cal A}
\Omega^1 \rightarrow \Omega^1 \otimes_{\cal A} \Gamma$ such that
\be
    \Psi_\nabla (\gamma \otimes_{\cal A} \d f)
  = \nabla (\gamma f) - (\nabla \gamma) f
    \qquad  \forall f \in {\cal A}, \, \gamma \in \Gamma   \; .
                              \label{Psi_nabla}
\ee
If we could turn this into a definition, we had a universal solution
to the problem we started with, the extension of connections on
two bimodules to the tensor product (over $\cal A$).
As a consequence of our assumptions, $\Psi_\nabla$ is defined on
$\Gamma \otimes_{\cal A} \Omega^1$, but the rhs of
(\ref{Psi_nabla}) must not respect that. This means that, depending on
the chosen differential calculus on $\cal A$, (\ref{Psi_nabla}) is only
well-defined for a special class of connections on $\Gamma$.
\vskip.2cm
\noindent
{\em Definition.} A connection $\nabla \, : \, \Gamma \rightarrow
\Omega^1 \otimes_{\cal A} \Gamma$ is called
{\em extensible}\footnote{In the notation of \cite{DV+Mass95} an
extensible connection is a `bimodule connection'.}
if it defines a map $\Psi_\nabla$ via (\ref{Psi_nabla}).
                                       \hfill  {\Large $\Box$}
\vskip.2cm
\noindent

Restrictions arise as follows. A relation $\sum h_k \, \d f_k=0$
in $\Omega^1$ implies $\sum \lbrack \nabla( \gamma \, h_k f_k)
- \nabla (\gamma h_k) \, f_k \rbrack =0$. If $\Omega^1$ is the
space of 1-forms of the universal first order differential calculus
on $\cal A$, there are no relations of the form $\sum h_k \, \d f_k=0$
and therefore every connection is extensible.
\vskip.2cm
\noindent
{\em Lemma \ref{Ext_con}.1.} For an extensible connection,
$\Psi_\nabla$ is an $\cal A$-bimodule homomorphism.
\vskip.1cm
\noindent
{\em Proof:}
\begin{eqnarray*}
     \Psi_\nabla (h \gamma \otimes_{\cal A} \d f)
 &=& \nabla ( h \gamma f) - \nabla (h \gamma) \, f  \\
 &=& \d h \otimes_{\cal A} (\gamma f) + h \, \nabla (\gamma f)
     - (\d h \otimes_{\cal A} \gamma) \, f - (h \nabla \gamma) \, f \\
 &=& h \, \lbrack \nabla (\gamma f) - (\nabla \gamma) \, f \rbrack \\
 &=& h \, \Psi_\nabla ( \gamma \otimes_{\cal A} \d f)         \\
     \Psi_\nabla (\gamma \otimes_{\cal A} (\d f) h)
 &=& \Psi_\nabla (\gamma \otimes_{\cal A} \d (f h) ) - \Psi_\nabla
     (\gamma \otimes_{\cal A} f \, \d h)    \\
 &=& \lbrack \nabla (\gamma f) - (\nabla \gamma) \, f \rbrack \, h \\
 &=& \Psi_\nabla ( \gamma \otimes_{\cal A} \d f) \, h
\end{eqnarray*}
for all $f,h \in {\cal A}$.                 \hfill    {\Large $\Box$}
\vskip.2cm
\noindent
What we have shown so far is summarized in the following proposition.
\vskip.2cm
\noindent
{\em Proposition \ref{Ext_con}.1.} For left $\cal A$-module connections
$\nabla, \nabla'$ on $\cal A$-bimodules $\Gamma, \Gamma'$ (with respect
to a first order differential calculus on $\cal A$) there exists a
connection on $\Gamma \otimes_{\cal A} \Gamma'$ of the form
\be
   \nabla = \Phi \circ (\nabla \otimes id_{\Gamma'})
   + \Psi \circ (id_{\Gamma} \otimes \nabla' )
\ee
with left $\cal A$-module homomorphisms $\Phi$ and $\Psi$ if and only
if $\nabla$ is extensible. The connection is then unique and given by
\be                \label{prod-conn}
   \nabla_\otimes := \nabla \otimes id_{\Gamma'}
 + (\Psi_\nabla \otimes id_{\Gamma'}) \circ (id_\Gamma \otimes \nabla')
\ee
with the $\cal A$-bimodule homomorphism $\Psi_\nabla$ defined via
(\ref{Psi_nabla}).                        \hfill    {\Large $\Box$}
\vskip.2cm
\noindent
If $\nabla \, : \, \Gamma \rightarrow \Omega^1
\otimes_{\cal A} \Gamma$ is an extensible connection on an
$\cal A$-bimodule $\Gamma$, a connection on the $n$-fold
tensor product $\Gamma^n$ of $\Gamma$ (over $\cal A$) is inductively
defined via
\be
   \nabla_{\otimes^n} := \nabla \otimes id_{\Gamma^{n-1}}
   + (\Psi_\nabla \otimes id_{\Gamma^{n-1}}) \circ (id_\Gamma
   \otimes \nabla^{(n-1)}) \; .
\ee
\vskip.2cm
\noindent
For an extensible connection we simply regard (\ref{Psi_nabla}) as the
definition of $\Psi_\nabla$. If, however, we choose some bimodule
homomorphism for $\Psi_\nabla$ on the lhs of (\ref{Psi_nabla}),
then this imposes further constraints on the connection. Corresponding
examples appeared in \cite{Mour95,DVMMM}.
\vskip.3cm
\noindent
{\em Proposition \ref{Ext_con}.2.} Let $\cal A$ be an associative
algebra and $(\Omega^1,\d )$ a first order differential calculus on
$\cal A$ which is inner, i.e., there is a 1-form $\rho$ such that
$\d f = \lbrack \rho , f \rbrack$ ($\forall f \in {\cal A}$). Let
$\sigma \, : \, \Omega^1 \otimes_{\cal A} \Omega^1 \rightarrow \Omega^1
\otimes_{\cal A} \Omega^1$ be a bimodule homomorphism.\footnote{A
possible choice is $\sigma \equiv 0$.}  \\
A linear connection is then extensible if and only if there exist
bimodule homomorphisms
\be
     V \, : \, \Omega^1 \otimes_{\cal A} \Omega^1 \rightarrow
            \Omega^1 \otimes_{\cal A} \Omega^1 \; , \qquad
     W \, : \, \Omega^1 \rightarrow \Omega^1 \otimes_{\cal A}
            \Omega^1
\ee
such that
\be                               \label{lc_decomp}
  \nabla \varphi = \nabla^\sigma \varphi
                   + V(\varphi \otimes_{\cal A} \rho) + W(\varphi)
\ee
where $\nabla^\sigma$ denotes the linear connection
associated\footnote{See also \cite{GMMM95}.}
with $\sigma$, defined by $\nabla^\sigma \varphi := \rho
\otimes_{\cal A} \varphi - \sigma(\varphi \otimes_{\cal A} \rho)$.
\vskip.1cm
\noindent
{\em Proof:}  \\
``$\Rightarrow$'': For an extensible $\nabla$ there is a bimodule
homomorphism $\Psi_\nabla \, : \, \Omega^1 \otimes_{\cal A} \Omega^1
\rightarrow \Omega^1 \otimes_{\cal A} \Omega^1$ such that
$\nabla (\varphi f) = (\nabla \varphi) f + \Psi_\nabla (\varphi
\otimes_{\cal A} \d f)$. Then
\begin{eqnarray*}
  \nabla^{\Psi_\nabla} \varphi := \rho \otimes_{\cal A} \varphi
            - \Psi_\nabla (\varphi \otimes_{\cal A} \rho)
\end{eqnarray*}
defines a linear connection and the difference $W := \nabla -
\nabla^{\Psi_\nabla}$ is a bimodule homomorphism. With the bimodule
homomorphism $V := \sigma - \Psi_\nabla$ we obtain the decomposition
(\ref{lc_decomp}).   \\
``$\Leftarrow$'': Assuming that (\ref{lc_decomp}) holds, we get
\begin{eqnarray*}
      \nabla (\varphi f) - (\nabla \varphi) \, f
  &=& \nabla^\sigma (\varphi f) - (\nabla^\sigma \varphi) f
      + V(\varphi f \otimes_{\cal A} \rho) - V(\varphi \otimes_{\cal A}
      \rho) \, f  \\
  &=& (\sigma - V)(\varphi \otimes_{\cal A} \d f)  \; .
\end{eqnarray*}
Since $\Psi_\nabla := \sigma - V$ is a bimodule homomorphism, $\nabla$
is extensible.                                \hfill  {\Large $\Box$}
\vskip.3cm

We should stress the following. The notion of extensibility of a
connection is based on the {\em ansatz} (\ref{ext-ansatz}). We cannot
exclude yet that there is a (more complicated) receipe to extend
connections to the tensor product of the modules on which they live,
without imposing restrictions on the connections. We have tried out
several modifications of (\ref{ext-ansatz}) without success.
\vskip.2cm

In \cite{Cunt+Quil95} a class of left $\cal A$-module connections on a
bimodule has been considered with additional right
$\cal A$-linearity.\footnote{The authors of \cite{Cunt+Quil95} call
such connections {\em left connections}. Furthermore, a
`connection on a bimodule' is defined in \cite{Cunt+Quil95} as a pair
of left and right connections.}
This is a subclass of extensible connections. For an algebra $\cal A$
with an inner first order differential calculus, there always exists
one particular connection of this kind, the {\em canonical (left
$\cal A$-module) connection} which is given by
$\nabla \varphi = \rho \otimes_{\cal A} \varphi$.
A complete characterization of such connections is obtained
in the following proposition.
\vskip.2cm
\noindent
{\em Proposition \ref{Ext_con}.3.} Let $\cal A$ be an associative
algebra and $(\Omega^1,\d )$ a first order differential calculus on
$\cal A$ which is inner (with a 1-form $\rho$). Then every left
$\cal A$-module connection which is also a right $\cal A$-module
homomorphism has the form
\be
   \nabla \varphi = \rho \otimes_{\cal A} \varphi + W(\varphi)
   \qquad \forall \varphi \in \Omega^1     \label{CQ-conn}
\ee
with a bimodule homomorphism $W$.
\vskip.1cm
\noindent
{\em Proof:} A left $\cal A$-module connection with the right
$\cal A$-module homomorphism property $\nabla(\varphi f) = (\nabla
\varphi) f$ is a special case of an extensible connection (with
$\Psi_\nabla = 0$). Proposition \ref{Ext_con}.2 then tells us that
$\nabla \varphi = \rho \otimes_{\cal A} \varphi + V(\varphi
\otimes_{\cal A} \rho) + W(\varphi)$ with bimodule homomorphisms
$V$ and $W$. The difference of two left $\cal A$-module connections
with right $\cal A$-linearity must be an $\cal A$-bimodule homomorphism
$\Omega^1 \rightarrow \Omega^1 \otimes_{\cal A} \Omega^1$. A simple
calculation using $\d f = \lbrack \rho , f \rbrack$ then shows that
$V$ has to vanish.                        \hfill  {\Large $\Box$}
\vskip.2cm
\noindent
The constraint imposed by (\ref{CQ-conn}) on a connection is very
restrictive. In section \ref{NgS3}.1 we have an example where the
canonical left $\cal A$-module connection turns out to be the only
left $\cal A$-module connection with right $\cal A$-linearity.
Further examples are provided by bicovariant first order differential
calculi on the quantum groups $GL_q(n)$.\footnote{A classification
of bicovariant differential calculi on the quantum general linear
groups $GL_q(n)$ has been obtained in \cite{Schm94}.}
In this case it has been shown \cite{GMMM95} that there is no
non-vanishing bimodule homomorphism $\Omega^1 \rightarrow \Omega^1
\otimes_{\cal A} \Omega^1$. All these calculi are inner \cite{Woro89}
(with a 1-form $\rho$). Hence, the canonical connection $\nabla \varphi
= \rho \otimes_{\cal A} \varphi$ is the only left $\cal A$-module
connection with right
$\cal A$-linearity according to Proposition \ref{Ext_con}.3.

\section{Connections and their duals}
\label{C_dual}
\setcounter{equation}{0}
Let $\cal A$ be an associative algebra and $\Gamma$ an $\cal
A$-bimodule. There are two natural ways to define a {\em dual}
of $\Gamma$, depending on whether its elements act from the left or
from the right on elements of $\Gamma$. Here we make the latter choice
(see the remark at the end of this section).
For the duality contraction $\langle \gamma , \mu \rangle$ where
$\gamma \in \Gamma$ and $\mu \in \Gamma^\ast$ (the dual of $\Gamma$),
we then have
\be                    \label{contr-lin}
  \langle f \gamma , \mu \rangle = f \, \langle \gamma , \mu \rangle
  \; , \qquad \langle \gamma , \mu f \rangle = \langle \gamma ,
  \mu \rangle \, f \; , \qquad
  \langle \gamma f , \mu \rangle = \langle \gamma , f \mu \rangle
\ee
for all $f \in {\cal A}$. For a left $\cal A$-module connection on
$\Gamma$ (with respect to some first order differential calculus
on $\cal A$ with space of 1-forms $\Omega^1$) its dual is a map
$\nabla^\ast \, : \, \Gamma^\ast \rightarrow \Gamma^\ast
\otimes_{\cal A} \Omega^1$ defined by
\be
   \langle \gamma , \nabla^\ast \mu \rangle := \d \langle \gamma ,
   \mu \rangle - \langle \nabla \gamma , \mu \rangle
\ee
where $\langle \gamma , \mu \otimes_{\cal A} \varphi \rangle := \langle
\gamma , \mu \rangle \, \varphi$ and $\langle \varphi \otimes_{\cal A}
\gamma , \mu \rangle := \varphi \, \langle \gamma , \mu \rangle$. With
these definitions we obtain
\be
     \langle \gamma , \nabla^\ast (\mu f) \rangle
 &=& \d \langle \gamma , \mu f \rangle -\langle \nabla \gamma , \mu f
     \rangle  \nonumber  \\
 &=& (\d \langle \gamma , \mu \rangle) \, f + \langle \gamma , \mu
     \rangle \, \d f - \langle \nabla \gamma , \mu \rangle \, f
     \nonumber  \\
 &=& \langle \gamma , \nabla^\ast \mu \rangle \, f + \langle \gamma ,
     \mu \otimes_{\cal A} \d f \rangle
\ee
and therefore
\be
  \nabla^\ast (\mu f) = (\nabla^\ast \mu) \, f
                       + \mu \otimes_{\cal A} \d f
\ee
which shows that $\nabla^\ast$ is a right $\cal A$-module connection.
For an extensible connection (see Appendix A) we now have the following
result.
\vskip.3cm
\noindent
{\em Proposition \ref{C_dual}.1.} If $\nabla \, : \, \Gamma \rightarrow
\Omega^1 \otimes_{\cal A} \Gamma$ is an extensible connection
with bimodule homomorphism $\Psi_\nabla$, then $\nabla^\ast$ is an
extensible connection with bimodule homomorphism $\Psi_{\nabla^\ast} \,
: \, \Omega^1 \otimes_{\cal A} \Gamma^\ast \rightarrow
\Gamma^\ast \otimes_{\cal A} \Omega^1$ given by
\be                         \label{dual-Psi}
  \langle \gamma , \Psi_{\nabla^\ast} (\varphi \otimes_{\cal A} \mu)
  \rangle := \langle \Psi_\nabla (\gamma \otimes_{\cal A} \varphi) ,
  \mu \rangle  \; .
\ee
\vskip.1cm
\noindent
{\em Proof:} It is easily checked that $\Psi_{\nabla^\ast}$ is
well-defined via (\ref{dual-Psi}) and that it is a bimodule
homomorphism. We still have to verify that $\Psi_{\nabla^\ast}$
satisfies the counterpart of (\ref{Psi_nabla}) for a right
$\cal A$-module connection,
\begin{eqnarray*}
  \langle \gamma , \nabla^\ast (f \mu) - f \, \nabla^\ast \mu \rangle
 &=& \langle \gamma , \nabla^\ast (f \mu) \rangle -
     \langle \gamma f , \nabla^\ast \mu \rangle \\
 &=& \d \langle \gamma , f \mu \rangle - \langle \nabla \gamma ,
     f \, \mu \rangle - \d \langle \gamma f , \mu \rangle
     + \langle \nabla (\gamma f) , \mu \rangle \\
 &=& \langle \nabla (\gamma f) - (\nabla \gamma) f , \mu \rangle
  = \langle \Psi_\nabla (\gamma \otimes_{\cal A} \d f) , \mu \rangle \;.
\end{eqnarray*}
\hspace*{1cm}                                    \hfill {\Large $\Box$}
\vskip.3cm

Let $\Gamma'$ be a left $\cal A$-module and $\nabla' \, : \, \Gamma'
\rightarrow \Omega^1 \otimes_{\cal A} \Gamma'$ a connection
on it. Its dual $\Gamma'^\ast$ is a right $\cal A$-module and
$\nabla'^\ast \, : \, \Gamma'^\ast \rightarrow \Gamma'^\ast
\otimes_{\cal A} \Omega^1$ defined as above is a connection
on $\Gamma'^\ast$.
In case we have on $\Gamma$ an extensible left $\cal A$-module
connection with a bimodule homomorphism $\Psi_\nabla$, we can define a
connection $\nabla_\otimes$ on
$\Gamma \otimes_{\cal A} \Gamma'$ in terms of the connections
on $\Gamma$ and $\Gamma'$ (see Appendix A). In the following we have
to assume that the dual module of $\Gamma \otimes_{\cal A} \Gamma'$ is
isomorphic to $\Gamma'^\ast \otimes_{\cal A} \Gamma^\ast$. This holds
in particular for modules of finite rank.
The duality contraction is then given by
\be
  \langle \gamma \otimes_{\cal A} \gamma' , \nu \otimes_{\cal A} \mu
  \rangle := \langle \gamma \, \langle \gamma' , \nu \rangle , \mu
  \rangle \; .
\ee
Now we have two different ways to define a connection on $\Gamma'^\ast
\otimes_{\cal A} \Gamma^\ast$, either as the dual of $\nabla_\otimes$,
i.e.,
\be
  \langle \gamma \otimes_{\cal A} \gamma' , (\nabla_\otimes)^\ast
  (\nu \otimes_{\cal A} \mu) \rangle := \d \langle \gamma
  \otimes_{\cal A} \gamma' , \nu \otimes_{\cal A} \mu \rangle -
  \langle \nabla_\otimes (\gamma \otimes_{\cal A} \gamma') , \nu
  \otimes_{\cal A} \mu \rangle \; ,
\ee
or as the `tensor product' of the duals $\nabla^\ast$ and
$\nabla'^\ast$, i.e.,
\be
 (\nabla^\ast)_\otimes (\nu \otimes_{\cal A} \mu) := (id_{\Gamma'}
 \otimes \Psi_{\nabla^\ast})(\nabla'^\ast \nu \otimes_{\cal A}
 \mu) + \nu \otimes_{\cal A} \nabla^\ast \mu \;.
\ee
Fortunately, both procedures lead to the same connection on
$\Gamma'^\ast \otimes_{\cal A} \Gamma^\ast$.
\vskip.3cm
\noindent
{\em Proposition \ref{C_dual}.2.}
\be
  (\nabla_\otimes)^\ast = (\nabla^\ast)_\otimes =: \nabla_\otimes^\ast
  \; .
\ee
{\em Proof:} Using
\begin{eqnarray*}
  \langle \gamma \otimes_{\cal A} \gamma' , (id_{\Gamma'}
  \otimes \Psi_{\nabla}^\ast)(\nabla'^\ast \nu \otimes_{\cal A}
  \mu) \rangle &=& \langle \Psi_\nabla (\gamma \otimes_{\cal A}
  \langle \gamma' , \nabla'^\ast \nu \rangle ) , \mu \rangle
\end{eqnarray*}
and
\begin{eqnarray*}
  \langle \Psi(\gamma \otimes_{\cal A} \langle \nabla' \gamma' , \nu
  \rangle) , \mu \rangle
  &=& \langle (\Psi \otimes id_{\Gamma'})
  (\gamma \otimes_{\cal A} \nabla' \gamma') , \nu \otimes_{\cal A}
  \mu \rangle \; ,
\end{eqnarray*}
a direct calculation shows that
\begin{eqnarray*}
 \langle \gamma \otimes_{\cal A} \gamma' , (\nabla_\otimes)^\ast(\nu
 \otimes_{\cal A} \mu) - (\nabla^\ast)_\otimes(\nu \otimes_{\cal A} \mu)
 \rangle =0   \; .
\end{eqnarray*}
\hspace*{1cm}                                    \hfill {\Large $\Box$}
\vskip.3cm
\noindent
{\em Remark.} Our choice among the two possible duals of $\Gamma$ is
related to our use of {\em left} $\cal A$-module connections. Let us
consider the alternative, the left dual $\Gamma'$ with contraction
$\langle \mu' , \gamma \rangle'$.\footnote{It still has to be clarified
whether the two duals, $\Gamma^\ast$ and $\Gamma'$ are isomorphic in
some (natural) sense.}
Then, in the expression $\langle \mu' , \nabla \gamma \rangle'$
the 1-form factor of $\nabla \gamma$ cannot be pulled out of the
contraction so that there is no (natural) way to define a dual of
a left $\cal A$-module connection. There is an exception, however. In
the special case of a linear connection, where $\Gamma = \Omega^1$, we
may indeed define a dual connection $\nabla'$ on the dual space
${\cal X}'$ of $\Omega^1$ via $\langle \nabla' X' , \varphi \rangle' =
\d \langle X' , \varphi \rangle' - \langle X' , \nabla \varphi
\rangle'$. Then $\nabla'$ is a left $\cal A$-module connection.
We emphasized earlier that, for a linear connection $\Omega^1
\rightarrow \Omega^1 \otimes_{\cal A} \Omega^1$, the two
$\Omega^1$-factors of the target space play very different roles.
So $\nabla'$ should only be taken seriously if there is a good
reason to forget about this fact.
Of course, if we consider {\em right} instead of left $\cal A$-module
connections, the correct contraction should be the primed one.
                                    \hfill    {\Large $\Box$}

\section{Coactions and extensions of invariant connections}
\label{Coact}
\setcounter{equation}{0}
Let $\cal A$ be a Hopf algebra with unit $\idty$ and coproduct $\Delta$,
$\Gamma$ an $\cal A$-bimodule and $\Gamma'$ a left $\cal A$-module.
$\Gamma$ and $\Gamma'$ are also assumed to be left $\cal A$-comodules
with coactions\footnote{In the following it will be sufficient to
consider a left coaction as a linear map $\Gamma \rightarrow {\cal A}
\otimes \Gamma$ such that $\Delta_{\Gamma}(f \gamma) = \Delta(f) \,
\Delta_{\Gamma}(\gamma)$. We will need a refinement in Proposition
\ref{Coact}.4 below, see the next footnote.}
\be
 \begin{array}{l@{\, : \;}l@{\; \rightarrow \;}l@{\qquad}l@{\, = \,}l}
   \Delta_\Gamma & \Gamma & {\cal A} \otimes \Gamma & \Delta_\Gamma
   (\gamma) & \sum_k f_k \otimes \gamma_k \\
   \Delta_{\Gamma'} & \Gamma' & {\cal A} \otimes \Gamma' &
   \Delta_{\Gamma'} (\gamma') & \sum_l f'_l \otimes \gamma'_l   \; .
 \end{array}
\ee
A left coaction on the tensor product $\Gamma \otimes_{\cal A}
\Gamma'$ is then given by
\be
   \Delta_{\Gamma \otimes_{\cal A} \Gamma'}(\gamma \otimes_{\cal A}
   \gamma') = \sum_{k,l} f_k \, f'_l \, \otimes \gamma_k
   \otimes_{\cal A} \gamma'_l
\ee
(see \cite{Woro89}, for example).
In the frequently used Sweedler notation \cite{Swee69}, this reads
\be
   \Delta_{\Gamma \otimes_{\cal A} \Gamma'}(\gamma \otimes_{\cal A}
   \gamma') = \gamma_{(-1)} \, \gamma'_{(-1)} \, \otimes \gamma_{(0)}
   \otimes_{\cal A} \gamma'_{(0)}
\ee
where $\Delta_\Gamma(\gamma) = \gamma_{(-1)} \otimes \gamma_{(0)}$.
\vskip.2cm

Let $(\Omega^1,\d )$ be a left-covariant first order differential
calculus on $\cal A$ and $\nabla \, : \, \Gamma \rightarrow \Omega^1
\otimes_{\cal A} \Gamma$ a left $\cal A$-module connection.
$\nabla$ is called {\em left-invariant} if
\be
   \Delta_{\Omega^1 \otimes_{\cal A} \Gamma} \circ \nabla
   = (id \otimes \nabla) \circ \Delta_{\Gamma}
\ee
where $\Delta_{\Omega^1 \otimes_{\cal A} \Gamma}$ is the left coaction
on $\Omega^1 \otimes_{\cal A} \Gamma$ induced by the left coactions on
$\Omega^1$ and $\Gamma$. As a consequence of this definition, if
$\gamma \in \Gamma$ is left-invariant (i.e., $\Delta_\Gamma (\gamma)
= \idty \otimes \gamma$) and if also $\nabla$ is left-invariant, then
$\nabla \gamma$ is left-invariant, i.e., $\Delta_{\Omega^1
\otimes_{\cal A} \Gamma} \nabla \gamma = \idty \otimes \nabla \gamma$.
\vskip.3cm
\noindent
{\em Proposition \ref{Coact}.1.} Let $\cal A$ be a Hopf algebra,
$\Omega^1$ a left-covariant differential calculus on $\cal A$,
$\Gamma$ an $\cal A$-bimodule and left $\cal A$-comodule with a
left-invariant extensible connection $\nabla$. Then the associated
bimodule homomorphism $\Psi_\nabla \, : \, \Gamma \otimes_{\cal A}
\Omega^1 \rightarrow \Omega^1 \otimes_{\cal A} \Gamma$ is also
left-invariant, i.e.,
\be
   \Delta_{\Omega^1 \otimes_{\cal A} \Gamma} \circ \Psi_\nabla
  = (id \otimes \Psi_\nabla) \circ \Delta_{\Gamma \otimes_{\cal A}
    \Omega^1}    \; .
\ee
\vskip.1cm
\noindent
{\em Proof:} Since all maps are $\Cx$-linear, it suffices to check
the invariance condition on elements of the form $\gamma
\otimes_{\cal A} \d f$ with $\gamma \in \Gamma$ and $f \in {\cal A}$,
\begin{eqnarray*}
 \Delta_{\Omega^1 \otimes_{\cal A} \Gamma} \circ \Psi_\nabla
     (\gamma \otimes_{\cal A} \d f)
 &=& \Delta_{\Omega^1 \otimes_{\cal A} \Gamma} (\nabla(\gamma f)
     - (\nabla \gamma) f)                 \\
 &=& (id \otimes \nabla) \circ \Delta_\Gamma (\gamma f)
     - ((id \otimes \nabla) \circ \Delta_\Gamma (\gamma)) \,
     \Delta(f)   \\
 &=& (id \otimes \nabla) (\Delta_\Gamma (\gamma) \, \Delta(f))
     - ((id \otimes \nabla) \circ \Delta_\Gamma (\gamma)) \,
     \Delta(f)   \\
 &=& \gamma_{(-1)} \, f_{(1)} \otimes \lbrack \nabla (\gamma_{(0)}
     \, f_{(2)}) - \nabla(\gamma_{(0)}) \, f_{(2)} \rbrack   \\
 &=& (id \otimes \Psi_\nabla) \circ (\gamma_{(-1)} \, f_{(1)}
     \otimes \gamma_{(0)} \otimes_{\cal A} \d f_{(2)})  \\
 &=& (id \otimes \Psi_\nabla) \circ \Delta_{\Gamma \otimes_{\cal A}
      \Omega^1} (\gamma \otimes_{\cal A} \d f)  \; .
\end{eqnarray*}
\hspace*{1cm}                           \hfill   {\Large  $\Box$}
\vskip.3cm
\noindent
{\em Proposition \ref{Coact}.2.} Let $\cal A$ be a Hopf algebra,
$\Omega^1$ a left-covariant differential calculus on $\cal A$, and
$\Gamma, \Gamma'$ two $\cal A$-bimodules which are also left
$\cal A$-comodules. Let $\nabla, \nabla'$ be left-invariant connections
on $\Gamma$ and $\Gamma'$, respectively. If $\nabla$ is extensible
(with associated bimodule homomorphism $\Psi_\nabla$), then
the product connection $\nabla_\otimes$ given by (\ref{prod-conn})
is a left-invariant connection on $\Gamma \otimes_{\cal A} \Gamma'$.
\vskip.1cm
\noindent
{\em Proof:}
\begin{eqnarray*}
 & & \Delta_{\Omega^1 \otimes_{\cal A} \Gamma \otimes_{\cal A}
     \Gamma'} \circ \nabla_\otimes (\gamma \otimes_{\cal A} \gamma') \\
 &=& \Delta_{(\Omega^1 \otimes_{\cal A} \Gamma) \otimes_{\cal A}
     \Gamma'} (\nabla \gamma \otimes_{\cal A} \gamma') +
     \Delta_{(\Omega^1 \otimes_{\cal A} \Gamma) \otimes_{\cal A}
     \Gamma'} \circ (\Psi_\nabla \otimes id) (\gamma \otimes_{\cal A}
     \nabla' \gamma')  \\
 &=& (id \otimes \nabla \otimes id) \circ \Delta_{\Gamma \otimes_{\cal A}
      \Gamma'} (\gamma \otimes_{\cal A} \gamma') +
     (id \otimes \Psi_\nabla \otimes id) \circ \Delta_{\Gamma
     \otimes_{\cal A} (\Omega^1 \otimes_{\cal A} \Gamma')}
     (\gamma \otimes_{\cal A} \nabla' \gamma')        \\
 &=& (id \otimes \nabla \otimes id) \circ \Delta_{\Gamma \otimes_{\cal A}
      \Gamma'} (\gamma \otimes_{\cal A} \gamma') +
     (id \otimes \Psi_\nabla \otimes id) \circ (id \otimes id
     \otimes \nabla') \circ \Delta_{\Gamma \otimes_{\cal A} \Gamma'}
     (\gamma \otimes_{\cal A} \gamma')        \\
 &=& (id \otimes \nabla_\otimes ) \circ \Delta_{\Gamma \otimes_{\cal
     A} \Gamma'} (\gamma \otimes_{\cal A} \gamma')  \; .
\end{eqnarray*}
\hspace*{1cm}                           \hfill   {\Large  $\Box$}
\vskip.3cm
\noindent
{\em Proposition \ref{Coact}.3.} Let $\cal A$ be a Hopf algebra and
$(\Omega^1, \d)$ a left-covariant (first order) differential calculus
on $\cal A$ which is inner, i.e., there is a 1-form $\rho$ such that
$\d f = \lbrack \rho , f \rbrack$ for all $f \in {\cal A}$. Let $\rho$
be left-invariant and $\Psi \, : \, \Omega^1 \otimes_{\cal A} \Omega^1
\rightarrow \Omega^1 \otimes_{\cal A} \Omega^1$ an $\cal A$-bimodule
homomorphism so that
\be
   \nabla^\Psi \varphi = \rho \otimes_{\cal A} \varphi
   - \Psi (\varphi \otimes_{\cal A} \rho)
\ee
defines a linear connection. Then $\nabla^\Psi$ is left-invariant
if and only if $\Psi$ is left-invariant.
\vskip.1cm
\noindent
{\em Proof:} \\
``$\Rightarrow$'': The connection $\nabla^\Psi$ is extensible and we
have $\Psi = \Psi_\nabla$. Hence $\Psi$ is left-invariant according to
Proposition \ref{Coact}.2.  \\
``$\Leftarrow$'': If $\Psi$ is left-invariant, then also $\nabla^\Psi$
since
\begin{eqnarray*}
 \Delta_{\Omega^1 \otimes_{\cal A} \Omega^1} \circ \nabla^\Psi \varphi
 &=& \Delta_{\Omega^1 \otimes_{\cal A} \Omega^1} (\rho \otimes_{\cal A}
     \varphi) - (id \otimes \Psi) \circ \Delta_{\Omega^1
     \otimes_{\cal A} \Omega^1} (\varphi \otimes_{\cal A} \rho) \\
 &=& \varphi_{(-1)} \otimes \rho \otimes_{\cal A} \varphi_{(0)}
     - (id \otimes \Psi) (\varphi_{(-1)} \otimes
     \varphi_{(0)} \otimes_{\cal A} \rho)   \\
 &=& \varphi_{(-1)} \otimes \nabla^\Psi \varphi_{(0)}   \\
 &=& (id \otimes \nabla^\Psi) \circ \Delta_{\Omega^1} (\varphi) \; .
\end{eqnarray*}
\hspace*{1cm}                           \hfill   {\Large  $\Box$}
\vskip.3cm

Let us now consider two $\cal A$-bimodules $\Gamma, \Gamma'$ with
{\em right} coactions
\be
 \begin{array}{l@{\, : \;}l@{\; \rightarrow \;}l@{\qquad}l@{\, = \,}l}
   {}_\Gamma\Delta & \Gamma & \Gamma \otimes {\cal A} & {}_\Gamma\Delta
   (\gamma) & \sum_k \gamma_k \otimes f_k  \\
   {}_{\Gamma'}\Delta & \Gamma' & \Gamma' \otimes {\cal A} & {}_{\Gamma'}
   \Delta (\gamma') & \sum_l \gamma'_l \otimes f'_l   \; .
 \end{array}
\ee
Then
\be
  {}_{\Gamma \otimes_{\cal A} \Gamma'}\Delta (\gamma \otimes_{\cal A}
   \gamma')  = \sum_{k,l} \gamma_k \otimes_{\cal A} \gamma'_l
   \otimes f_k \, f'_l
\ee
defines a right coaction on $\Gamma \otimes_{\cal A} \Gamma'$.
\vskip.2cm

Let $(\Omega^1, \d)$ be a {\em right}-covariant first order
differential calculus on $\cal A$. If an $\cal A$-bimodule $\Gamma$ has
a right coaction $_{\Gamma}\Delta \, : \, \Gamma \rightarrow \Gamma
\otimes {\cal A}$, then {\em right-invariance} of a connection $\nabla$
is defined by
\be
  {}_{\Omega^1 \otimes_{\cal A} \Gamma}\Delta \circ \nabla
   = (\nabla \otimes id) \circ {}_{\Gamma}\Delta
\ee
where ${}_{\Omega^1 \otimes_{\cal A} \Gamma} \Delta$ is the right
coaction on $\Omega^1 \otimes_{\cal A} \Gamma$ induced by the right
coactions on $\Omega^1$ and $\Gamma$. With these notions, the last
three propositions in this section remain valid if `left' is everywhere
replaced by `right' (with the exception that we still consider left
$\cal A$-module connections).
\vskip.3cm

In the following we demonstrate that invariance properties of
connections are also transfered to their duals. First we establish
the existence of a left coaction on the dual of a bimodule with
a left coaction.
\vskip.3cm
\noindent
{\em Proposition \ref{Coact}.4.} Let $\cal A$ be a Hopf algebra
and $\Gamma$ a left-covariant bimodule over $\cal A$ with coaction
$\Delta_\Gamma$.\footnote{A left-covariant bimodule over a
Hopf algebra $\cal A$ is an $\cal A$-bimodule $\Gamma$ together
with a map (coaction) $\Delta_\Gamma \, : \, \Gamma \rightarrow
{\cal A} \otimes \Gamma$ such that $\Delta_\Gamma(f \gamma f')
= \Delta(f) \, \Delta_\Gamma(\gamma) \, \Delta(f')$ for all
$f,f' \in {\cal A}, \, \gamma \in \Gamma$. It has to satisfy the
equations $(\Delta \otimes id) \circ \Delta_\Gamma = (id \otimes
\Delta_\Gamma) \circ \Delta_\Gamma$ and $(\epsilon \otimes id)
\circ \Delta_\Gamma = id$ where $\epsilon$ is the counit.
See also \cite{Woro89}. }
Then the dual module $\Gamma^\ast$ has a unique left-covariant bimodule
structure with coaction $\Delta_{\Gamma^\ast} \, : \, \Gamma^\ast
\rightarrow {\cal A} \otimes \Gamma^\ast$ such that
\be                    \label{dual-lcoact}
 (id \otimes \langle \; , \; \rangle ) \circ \Delta_{\Gamma
 \otimes_{\cal A} \Gamma^\ast} = \Delta \circ \langle \; , \; \rangle
\ee
where $\langle \; , \; \rangle$ denotes the contraction mapping
$\Gamma \otimes_{\cal A} \Gamma^\ast \rightarrow {\cal A}$.
\vskip.1cm
\noindent
{\em Proof:} According to Theorem 2.1 in \cite{Woro89} a left-covariant
bimodule $\Gamma$ has a left $\cal A$-module basis of left-invariant
elements $\lbrace \gamma^i \rbrace$ (where $i$ runs through some
index set). Let $\lbrace \mu_j \rbrace$ be the dual basis of
$\Gamma^\ast$. Assuming the existence of $\Delta_{\Gamma^\ast}$,
(\ref{dual-lcoact}) leads to $\mu_{j(-1)} \otimes \langle \gamma^i ,
\mu_{j(0)} \rangle = \delta^i_j \, \idty \otimes \idty$
which implies $\mu_{j(-1)} \in \Cx$. Now $(\epsilon \otimes id) \circ
\Delta_\Gamma = id$ shows that the $\mu_j$ are left-invariant,
i.e., $\Delta_{\Gamma^\ast}(\mu_j) = \idty \otimes \mu_j$. As a
consequence, the coaction is unique.

Let us now prove the existence of the coaction on $\Gamma^\ast$.
According to Theorem 2.1 in \cite{Woro89} there are maps $F^i{}_j \,
: \, {\cal A} \rightarrow {\cal A} \,$ such that
\begin{eqnarray*}
  \gamma^i \, f = \sum_k F^i{}_k(f) \, \gamma^k  \; , \quad
  \Delta(F^i{}_k(f)) = (id \otimes F^i{}_k) \, \Delta(f) \qquad
  \forall f \in {\cal A}  \; .
\end{eqnarray*}
Then
\begin{eqnarray*}
     \langle \gamma^i , f \, \mu_j \rangle
  &=& \langle \gamma^i \, f , \mu_j \rangle
   = \sum_k F^i{}_k(f) \, \langle \gamma^k , \mu_j \rangle
   = F^i{}_j(f)   \\
  &=& \sum_k \langle \gamma^i , \mu_k \rangle \, F^k{}_j(f)
   = \langle \gamma^i , \sum_k \mu_k \, F^k{}_j(f) \rangle
\end{eqnarray*}
implies $f \, \mu_j = \sum_k \mu_k \, F^k{}_j(f)$.
Now we define the coaction $\Delta_{\Gamma^\ast}$ on the basis $\mu_j$
by $ \Delta_{\Gamma^\ast} (\mu_j) := \idty \otimes \mu_j $
and extend it via $\Delta_{\Gamma^\ast} (\mu_j \, f) :=
\Delta_{\Gamma^\ast} (\mu_j) \, \Delta(f)$ for all $f \in {\cal A}$.
Then
\begin{eqnarray*}
      \Delta_{\Gamma^\ast} (f \, \mu_j)
  &=& \Delta_{\Gamma^\ast} ( \sum_k \mu_k \, F^k{}_j(f) )
   =  \sum_k \Delta_{\Gamma^\ast} (\mu_k) \, \Delta(F^k{}_j(f))  \\
  &=& \sum_k (\idty \otimes \mu_k) \, (id \otimes F^k{}_j) \,
      \Delta(f)
   =  \sum_k f_{(1)} \otimes \mu_k \, F^k{}_j(f_{(2)})
   =  \sum_k f_{(1)} \otimes f_{(2)} \, \mu_j                     \\
  &=& \Delta(f) \, \Delta_\Gamma(\mu_j)       \; .
\end{eqnarray*}
It is now sufficient to verify the remaining defining properties of
a left-covariant bimodule on the left-invariant basis elements
$\mu_j$, and furthermore (\ref{dual-lcoact}) on $\lbrace \gamma^i
\rbrace$ and $\lbrace \mu_j \rbrace$. We leave this to the reader.
                                            \hfill  {\Large $\Box$}
\vskip.3cm

After some preparations in the following Lemma, we prove that
left-invariance of a connection on a left-covariant bimodule translates
to invariance of the dual connection which lives on the dual
left-covariant bimodule.
\vskip.3cm
\noindent
{\em Lemma \ref{Coact}.1.}
\be
   (id \otimes \langle \; , \; \rangle) \circ \Delta_{\Gamma
   \otimes_{\cal A} \Gamma^\ast \otimes_{\cal A} \Omega^1}
   (\gamma \otimes_{\cal A} \hat{\mu} ) &=&
   \Delta_{\Omega^1} \langle \gamma , \hat{\mu} \rangle
     \qquad \forall \hat{\mu} \in \Gamma^\ast \otimes_{\cal A}
            \Omega^1   \\
   (id \otimes \langle \; , \; \rangle) \circ \Delta_{\Omega^1
   \otimes_{\cal A} \Gamma \otimes_{\cal A} \Gamma^\ast}
   (\hat{\gamma} \otimes_{\cal A} \mu) &=&
   \Delta_{\Omega^1} \langle \hat{\gamma} , \mu \rangle
   \qquad \forall \hat{\gamma} \in \Omega^1 \otimes_{\cal A} \Gamma
             \; .
\ee
\vskip.1cm
\noindent
{\em Proof:} This is a straightforward calculation using
(\ref{dual-lcoact}).     \hfill  {\Large $\Box$}

\vskip.3cm
\noindent
{\em Proposition \ref{Coact}.5.} Let $\cal A$ be a Hopf algebra,
$\Gamma$ a left-covariant bimodule over $\cal A$ and $(\Omega^1, \d)$
a left-covariant first order differential calculus on $\cal A$.
If $\nabla \, : \, \Gamma \rightarrow \Omega^1 \otimes_{\cal A}
\Gamma$ is left-invariant, then the dual connection $\nabla^\ast
\, : \, \Gamma^\ast \rightarrow \Gamma^\ast \otimes_{\cal A}
\Omega^1$ is also left-invariant.
\vskip.1cm
\noindent
{\em Proof:} We have to show that
\begin{eqnarray*}
  \Delta_{\Gamma^\ast \otimes_{\cal A} \Omega^1} \circ \nabla^\ast
  = (id \otimes \nabla^\ast) \circ \Delta_{\Gamma^\ast}  \; .
\end{eqnarray*}
Let $\lbrace \gamma^i \rbrace$ be a left-invariant left $\cal
A$-module basis of $\Gamma$ \cite{Woro89}. We introduce mappings
$C^i \, : \, \Gamma^\ast \otimes_{\cal A} \Omega^1 \rightarrow
\Omega^1$, $\mu \otimes_{\cal A} \varphi \mapsto \langle \gamma^i ,
\mu \rangle \, \varphi$. Then
\begin{eqnarray*}
  & & (id \otimes C^i) (id \otimes \nabla^\ast) \, \Delta_{\Gamma^\ast}
      (\mu)
   =  \mu_{(-1)} \otimes \langle \gamma^i , \nabla^\ast \mu_{(0)}
      \rangle
   =  \mu_{(-1)} \otimes \lbrack \d \langle \gamma^i , \mu_{(0)} \rangle
      - \langle \nabla \gamma^i , \mu_{(0)} \rangle \rbrack   \\
  &=& \lbrack (id \otimes \d) \circ (id \otimes \langle \; , \;
      \rangle) - (id \otimes \langle \; , \; \rangle) \circ
      (id \otimes \nabla \otimes id) \rbrack \circ \Delta_{\Gamma
      \otimes_{\cal A} \Gamma^\ast} (\gamma^i \otimes_{\cal A} \mu) \\
  &=& (id \otimes \d) \circ \Delta \langle \gamma^i , \mu \rangle
      - (id \otimes \langle \; , \; \rangle) \circ
      \Delta_{\Omega^1 \otimes_{\cal A} \Gamma \otimes_{\cal A}
      \Gamma^\ast} \circ (\nabla \otimes id) (\gamma^i \otimes_{\cal A}
      \mu)           \\
  &=& \Delta_{\Omega^1} \, ( \d \langle \gamma^i , \mu \rangle
      -  \langle \nabla \gamma^i , \mu \rangle )
      \, = \, \Delta_{\Omega^1} \langle \gamma^i , \nabla^\ast \mu
      \rangle
   =  (id \otimes \langle \; , \; \rangle) \circ
      \Delta_{\Gamma \otimes_{\cal A} \Gamma^\ast \otimes_{\cal A}
      \Omega^1} (\gamma^i \otimes_{\cal A} \nabla^\ast \mu) \\
  &=& (\nabla^\ast \mu)_{(-1)} \otimes \langle \gamma^i ,
      (\nabla^\ast \mu)_{(0)} \rangle
   =  (id \otimes C^i) \circ \Delta_{\Gamma^\ast \otimes_{\cal A}
      \Omega^1} (\nabla^\ast \mu)
\end{eqnarray*}
using the left-invariance of $\nabla \gamma^i$ and Lemma
\ref{Coact}.1.
It remains to show that $(id \otimes C^i)(\xi) = 0$ for all $i$
implies $\xi = 0$ (where $\xi \in {\cal A} \otimes \Gamma^\ast
\otimes_{\cal A} \Omega^1$). $\xi$ has an expression
\begin{eqnarray*}
  \xi = \sum_{\alpha,j,r} f_\alpha \otimes \mu_j \otimes_{\cal A}
          \varphi_r \; \xi_{\alpha j r}
%      = \sum_{j,r} \lbrack \idty \otimes (\mu_j \otimes_{\cal A}
%        \varphi_r) \rbrack \, ( \sum_{\alpha} f_\alpha \otimes
%        \xi_{\alpha j r} )
\end{eqnarray*}
with $\xi_{\alpha j r} \in {\cal A}$. Here $\lbrace \mu_j \rbrace$
is the basis of $\Gamma^\ast$ dual to $\lbrace \gamma^i \rbrace$ and
$\lbrace \varphi_r \rbrace$ is a right $\cal A$-module basis of
$\Omega^1$. Evaluation of $(id \otimes C^i)(\xi) = 0$ now leads to
$\sum_\alpha f_\alpha \otimes \xi_{\alpha i r} = 0$ for all $i$ and
all $r$. Hence $\xi = 0$.          \hfill  {\Large $\Box$}

\vskip.3cm
\noindent
Corresponding results are obtained for a right-covariant bimodule
$\Gamma$ with coaction ${}_\Gamma \Delta$ and right-invariant
connections on it. In this case (\ref{dual-lcoact}) is replaced
by
\be
 (\langle \; , \; \rangle \otimes id) \circ {}_{\Gamma \otimes_{\cal A}
 \Gamma^\ast}\Delta = \Delta \circ \langle \; , \; \rangle   \; .
\ee

\section{Two-sided connections}
\setcounter{equation}{0}
The problem of extensibility of a connection discussed in Appendix
A disappears if we modify its definition as follows.
\vskip.2cm
\noindent
{\em Definition.} A {\em two-sided connection}\footnote{See also
\cite{DV+Mich95,DV+Mass95} for related structures.}
on an $\cal A$-bimodule $\Gamma$ is a map $\nabla \, : \, \Gamma
\rightarrow (\Omega^1 \otimes_{\cal A} \Gamma) \oplus (\Gamma
\otimes_{\cal A} \Omega^1)$ such that
\be
   \nabla (f \gamma f') = \d f \otimes_{\cal A} \gamma \, f'
   + f \, \gamma \otimes_{\cal A} \d f' + f \, (\nabla \gamma) \, f'
\ee
for all $f,f' \in {\cal A}$ and $\gamma \in \Gamma$.
\vskip.2cm
\noindent
The difference of two such connections is a bimodule homomorphism. The
following example demonstrates that the concept of a two-sided
connection is much more restrictive than that of the usual one.
\vskip.3cm
\noindent
{\em Example.} For a first order differential calculus which is
inner, i.e., there is a 1-form $\rho$ such that $\d f = \lbrack
\rho , f \rbrack$ for all $f \in {\cal A}$,
\be
   \nabla \varphi := \rho \otimes_{\cal A} \varphi - \varphi
                    \otimes_{\cal A} \rho
\ee
defines a two-sided linear connection. In the particular case of the
three-dimensional bicovariant differential calculus on ${\cal S}_3$,
we observed in section \ref{NgS3} that there is no nontrivial bimodule
homomorphism $\Omega^1 \rightarrow \Omega^1 \otimes_{\cal A} \Omega^1$.
Hence, the two-sided connection defined above is the only one in this
case.                                        \hfill   {\Large $\Box$}
\vskip.3cm

A two-sided connection extends to a map $\Omega \otimes_{\cal A} \Gamma
\otimes_{\cal A} \Omega \rightarrow \Omega \otimes_{\cal A} \Gamma
\otimes_{\cal A} \Omega$ via
\be
   \nabla (\varphi \gamma \varphi') = (\d \varphi) \, \gamma
   \, \varphi' + (-1)^r \, (\nabla \gamma) \, \varphi'
   + (-1)^{r+s} \, \varphi \, \gamma \, \d \varphi'
\ee
where $\varphi \in \Omega^r$ and $\gamma \in \bigoplus_{k=0}^s
\Omega^k \otimes_{\cal A} \Gamma \otimes_{\cal A} \Omega^{s-k}$.
The curvature of $\nabla$ then turns out to be an $\cal A$-bimodule
homomorphism, i.e.,
\be
    \nabla^2 (f \gamma f') = f \, (\nabla^2 \gamma) \, f'
    \qquad \forall f,f' \in {\cal A}, \, \gamma \in \Gamma \, ,
\ee
a nice property not shared, in general, by ordinary connections.

\section{Invariant tensor fields on a finite group}
\label{Invtensor}
\setcounter{equation}{0}
{}From two $\cal A$-bimodules $\Gamma$ and $\Gamma'$ we can build
the tensor product $\Gamma \otimes_{\cal A} \Gamma'$. If both
modules carry a (left or right) $\cal A$-comodule structure, there is
a comodule structure on the tensor product space (see Appendix C).
In case of left comodules, the left-invariance condition for a tensor
field $\alpha \in \Gamma \otimes_{\cal A} \Gamma'$ reads
$\Delta_{\Gamma \otimes_{\cal A} \Gamma'} (\alpha) = \idty \otimes
\alpha$. For right comodules this is replaced by the right-invariance
condition ${}_{\Gamma \otimes_{\cal A} \Gamma'}\Delta (\alpha)
= \alpha \otimes \idty$.
In the following we consider a bicovariant (first order) differential
calculus on a finite group $G$. Besides being an $\cal A$-bimodule, the
space $\Omega^1$ is then a left and right $\cal A$-comodule.
Each tensor field $\alpha \in \Omega^1 \otimes_{\cal A} \Omega^1$
can be written as
\be                                       \label{alpha}
   \alpha = \alpha_{\un{g},\un{g}'} \, \theta^{\un{g}}
            \otimes_{\cal A} \theta^{\un{g}'}
\ee
where summations run over the set $\hat{G} = \lbrace g \in G \mid
\theta^g \neq 0 \rbrace$.
Left-invariance of $\alpha$ then means $\alpha_{g,g'} \in \Cx$.
Bi-invariance leads to the additional condition
\be
   \alpha_{\un{g},\un{g}'} \, \theta^{ad(h) \un{g}} \otimes_{\cal
   A} \theta^{ad(h) \un{g}'} = \alpha_{\un{g},\un{g}'} \,
   \theta^{\un{g}} \otimes_{\cal A} \theta^{\un{g}'}
   \qquad \forall h \in G  \; .
\ee
For fixed $h \in G$, the map $ad(h) \, : \, G \rightarrow G$ is a
bijection. Hence
\be
   \alpha_{\un{g},\un{g}'} \, \theta^{ad(h) \un{g}} \otimes_{\cal
   A} \theta^{ad(h) \un{g}'} = \alpha_{ad(h^{-1}) \un{k},
   ad(h^{-1}) \un{k}'} \, \theta^{\un{k}} \otimes_{\cal A}
   \theta^{\un{k}'}
\ee
and the bi-invariance condition becomes
\be
   \alpha_{ad(h)g, ad(h)g'} = \alpha_{g,g'} \in \Cx
   \qquad  \forall g,g' \in \hat{G}, \, h \in G  \; .
\ee
For a bicovariant differential calculus with bimodule isomorphism
$\sigma$, the condition for a tensor field $\alpha$ to be s-symmetric
is
\be
  \alpha_{g,g'} = \alpha_{g', ad(g')g}  \qquad \forall g,g' \in
  \hat{G}  \; .
\ee
$\alpha$ is s-antisymmetric iff
\be
  \alpha_{g,g'} = - \alpha_{g', ad(g')g}  \qquad \forall g,g' \in
  \hat{G}  \; .
\ee
\vskip.3cm
\noindent
{\em Example.} Let us consider ${\cal S}_3$ with the universal (first
order) differential calculus (see section \ref{NgS3}).
In matrix notation, the coefficients of an s-symmetric tensor field
$\alpha$, as given by (\ref{alpha}), must have the form
\be
   \left( \alpha_{g,g'} \right) = \left( \begin{array}{ccccc}
   \alpha_1 & \alpha_4 & \alpha_5 & \beta_1 & \beta_2  \\
   \alpha_5 & \alpha_2 & \alpha_4 & \beta_3 & \beta_1  \\
   \alpha_4 & \alpha_5 & \alpha_3 & \beta_2 & \beta_3  \\
   \beta_2  & \beta_1  & \beta_3  & \gamma_1& \gamma_3 \\
   \beta_1  & \beta_3  & \beta_2  & \gamma_3& \gamma_2
   \end{array} \right)
\ee
where the entries are (arbitrary) elements of $\cal A$ (respectively
constants if $\alpha$ is left-invariant). Rows and columns are
arranged, respectively, according to the index sequence $\lbrace a,b,c,ab,
ba \rbrace$. For an s-antisymmetric tensor field we obtain
\be
   \left( \alpha_{g,g'} \right) = \left( \begin{array}{ccccc}
   0        & 0        & 0        & \beta_1 & \beta_2  \\
   0        & 0        & 0        & \beta_3 & \beta_1  \\
   0        & 0        & 0        & \beta_2 & \beta_3  \\
   -\beta_2  & -\beta_1  & -\beta_3  & 0    & \gamma \\
   -\beta_1  & -\beta_3  & -\beta_2  & -\gamma  & 0
   \end{array} \right)  \; .
\ee
For a w-symmetric tensor field we find
\be
   \left( \alpha_{g,g'} \right) = \left( \begin{array}{ccccc}
   \alpha_1 & \alpha_4 & \alpha_5 & \beta_1 & \beta_1'  \\
   \alpha_7 & \alpha_2 & \alpha_6 & \beta_1'' & \beta_3  \\
   \alpha_8 & \alpha_9 & \alpha_3 & \beta_3' & \beta_3''  \\
   \beta_1'-\beta_2'+\beta_3' & \beta_2  & \beta_2''  & \gamma_1&
                                                        \gamma_3 \\
   \beta_1-\beta_2+\beta_3  & \beta_1''-\beta_2''+\beta_3''
   & \beta_2'  & \gamma_3& \gamma_2
   \end{array} \right)
\ee
and the coefficients of a w-antisymmetric tensor field are given by
\be
   \left( \alpha_{g,g'} \right) = \left( \begin{array}{ccccc}
   0 & \alpha_1 & -\alpha_3-\alpha_4 & \beta_1 & \beta_1'  \\
   \alpha_3 & 0 & \alpha_2 & \beta_1'' & \beta_3  \\
   -\alpha_1-\alpha_2 & \alpha_4 & 0 & \beta_3' & \beta_3''  \\
   -\beta_1'-\beta_2'-\beta_3' & \beta_2  & \beta_2''  & 0 & \gamma \\
   -\beta_1-\beta_2-\beta_3  & -\beta_1''-\beta_2''-\beta_3''
   & \beta_2'  & -\gamma& 0
   \end{array} \right) \; .
\ee
$\alpha$ is bi-invariant iff the coefficient matrix has the form
\be
   \left( \alpha_{g,g'} \right) = \left( \begin{array}{ccccc}
   \alpha  & \alpha' & \alpha' & \beta  & \beta   \\
   \alpha' & \alpha  & \alpha' & \beta  & \beta   \\
   \alpha' & \alpha' & \alpha  & \beta  & \beta   \\
   \beta   & \beta   & \beta   & \gamma & \gamma' \\
   \beta   & \beta   & \beta   & \gamma'& \gamma
   \end{array} \right)
   +  \left( \begin{array}{ccccc}
   0        & 0        & 0        & \beta'  & \beta'   \\
   0        & 0        & 0        & \beta'  & \beta'   \\
   0        & 0        & 0        & \beta'  & \beta'   \\
   -\beta'  & -\beta'  & -\beta'  & 0    &  0          \\
   -\beta'  & -\beta'  & -\beta'  & 0    &  0
   \end{array} \right)
\ee
with constants $\alpha,\alpha',\beta,\beta',\gamma,\gamma'$.
As expressed above, it turns out to be a sum of s-symmetric and
s-antisymmetric tensors.

\section{Finite group actions on a finite set}
\setcounter{equation}{0}
Within the framework of noncommutative geometry of finite sets
one can also formulate the notion of covariance with
respect to a group action on a finite set. Let $M = \{x,y,\ldots\}$ be
this set and $G$ a finite group acting on $M$ from the left,
\be
   G \times M \rightarrow M \qquad (g,x) \mapsto g \cdot x  \; .
\ee
For $g,g'\in G$ and $x \in M$ we have $(gg') \cdot x = g \cdot (g'
\cdot x)$. The action of the neutral element $e \in G$ is trivial,
i.e., $e \cdot x = x$ for all $x \in M$.
We denote the algebra of $\Cx$-valued functions on $M$ and $G$ by
${\cal H}$ and ${\cal A}$, respectively. ${\cal A}$ is a Hopf algebra
over $\Cx$. The group action induces a
{\em left coaction} $\Delta_{\cal H}: {\cal H} \rightarrow {\cal A}
\otimes {\cal H}$ via
\begin{equation}
         \Delta_{\cal H}(f)(g,x) = f(g \cdot x)  \; .
\end{equation}
Since $\Delta_{\cal H}$ is compatible with the multiplication in
${\cal H}$, the latter is turned into a (left) ${\cal A}$-{\em comodule
algebra}. In particular,
\be
   \Delta_{\cal H}(e_x) =  \sum_{g \in G} e_g \otimes e_{g^{-1}
   \cdot x}  \; .
\ee
A (first order) differential calculus on $M$ (or $\cal H$) with space of
1-forms $\Omega^1$ is called $G$-{\em covariant}
iff there is a linear map $\ldel \, : \, \Omega^1 \rightarrow {\cal A}
\otimes \Omega^1$ such that
\be
  \ldel (f \, \varphi\, h) = \Delta_{\cal H}(f) \, \ldel (\varphi)
  \, \Delta_{\cal H}(h)   \qquad
  \forall f \in {\cal A}, \; h \in {\cal H}
\ee
and
\be
      \ldel \circ \d = (id \otimes \d) \circ \Delta_{\cal H} \; .
\ee
As a consequence,
\be
   \ldel (e_{x,y}) = \sum_{g \in G} e_{g^{-1}} \otimes
                     e_{g \cdot x, g \cdot y}  \; .
\ee
We obtain all $G$-covariant differential calculi on $M$ by deleting
sets of arrows from the universal graph (the digraph which corresponds
to the universal differential calculus on $M$). These correspond to
$G$-orbits in $(M \times M)'$. A (nontrivial) $G$-covariant differential
calculus is called {\em irreducible} if it belongs to a single orbit.
All (nontrivial) differential calculi are then obtained as unions
of irreducible ones.
\vskip.2cm
\noindent
{\em Example 1.} Let $M$ be a finite set with $n$ elements and
$G = {\cal S}_n$, the symmetric group. Obviously, the action of
${\cal S}_n$ on $(M \times M)'$ is transitive, i.e., all $(x,y)$
where $x \neq y$ belong to the same ${\cal S}_n$-orbit. Therefore,
the only ${\cal S}_n$-covariant (first order) differential calculi
on $M$ are the universal and the trivial one.
\vskip.2cm
\noindent
{\em Example 2.} Instead of the action of the whole symmetric group
we may consider actions of subgroups of ${\cal S}_n$. For $n=3$, for
example, we have the (nontrivial) subgroups
\begin{eqnarray}
   G_1 = \{e,a\} \, , \quad
   G_2 = \{e,b\} \, , \quad
   G_3 = \{e,c\} \, , \quad
   G_4 = \{e,ab,ba\} \; .
\end{eqnarray}
Denoting the points of $M$ by $1,2,3$, we can calculate the orbits
in $(M \times M)'$. For the action of $G_1$, we obtain
\begin{eqnarray}
   {\cal O}_1 = \{(1,2),(2,1)\} \, , \quad
   {\cal O}_2 = \{(1,3),(2,3)\} \, , \quad
   {\cal O}_3 = \{(3,1),(3,2)\} \; .
\end{eqnarray}
The graphs which determine the irreducible calculi are depicted in
Fig. 3.

\unitlength1.cm
\begin{picture}(12.,1.9)(-4.2,-0.3)
\thicklines
\put(0.,1.) {\circle*{0.1}}
\put(1.,0.) {\circle*{0.1}}
\put(-1.,0.) {\circle*{0.1}}
\put(0.8,0.1) {\vector(-1,1){0.7}}
\put(0.2,0.9) {\vector(1,-1){0.7}}
\put(-0.1,0.8) {\vector(-1,-1){0.7}}
\put(-0.9,0.2) {\vector(1,1){0.7}}
\put(-1.3,-0.2){$1$}
\put(-0.1,1.2){$3$}
\put(1.2,-0.2){$2$}
\put(3.,1.) {\circle*{0.1}}
\put(4.,0.) {\circle*{0.1}}
\put(2.,0.) {\circle*{0.1}}
\put(2.2,0.05) {\vector(1,0){1.6}}
\put(3.8,-0.07) {\vector(-1,0){1.6}}
\put(3.1,0.9) {\vector(1,-1){0.8}}
\put(2.9,0.9) {\vector(-1,-1){0.8}}
\put(6.,1.) {\circle*{0.1}}
\put(7.,0.) {\circle*{0.1}}
\put(5.,0.) {\circle*{0.1}}
\put(5.2,0.05) {\vector(1,0){1.6}}
\put(6.8,-0.07) {\vector(-1,0){1.6}}
\put(5.1,0.1) {\vector(1,1){0.8}}
\put(6.9,0.1) {\vector(-1,1){0.8}}
\end{picture}

\begin{center}
\begin{minipage}{10.5cm}
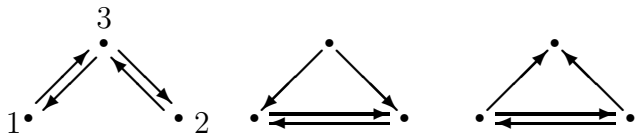
\centerline{\bf Fig. 3}
\vskip.1cm
\noindent
The digraphs coresponding to the irreducible $G_1$-covariant first
order differential calculi on a 3-point set.
\end{minipage}
\end{center}
\vskip.3cm
\noindent
In the case of $G_4$ acting on $M$ we obtain
\begin{eqnarray}
   {\cal O}_1 = \{(1,2),(2,3),(3,1)\} \, , \quad
   {\cal O}_2 = \{(1,3),(2,1),(3,2)\} \; .
\end{eqnarray}
The graphs corresponding to irreducible calculi are displayed in Fig. 4.

\unitlength1.cm
\begin{picture}(12.,1.8)(-2.7,-0.3)
\thicklines
\put(3.,1.) {\circle*{0.1}}
\put(4.,0.) {\circle*{0.1}}
\put(2.,0.) {\circle*{0.1}}
\put(3.9,0.) {\vector(-1,0){1.8}}
\put(3.1,0.9) {\vector(1,-1){0.8}}
\put(2.1,0.1) {\vector(1,1){0.8}}
\put(6.,1.) {\circle*{0.1}}
\put(7.,0.) {\circle*{0.1}}
\put(5.,0.) {\circle*{0.1}}
\put(5.1,0.) {\vector(1,0){1.8}}
\put(6.9,0.1) {\vector(-1,1){0.8}}
\put(5.9,0.9) {\vector(-1,-1){0.8}}
\end{picture}

\begin{center}
\begin{minipage}{10.5cm}
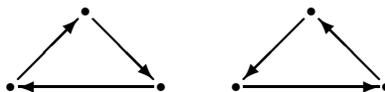
\centerline{\bf Fig. 4}
\vskip.1cm
\noindent
The digraphs corresponding to the irreducible $G_4$-covariant first
order differential calculi on a 3-point set.
\end{minipage}
\end{center}
\vskip.2cm

Of course, one can proceed with the formalism by defining invariance
of tensors and connections on $M$. All this will be explored in
detail in a separate work.

\end{appendix}

\end{document}